\documentclass[11pt]{amsart}
\usepackage{amsmath,amsthm,amsfonts,amssymb,stmaryrd,mathrsfs,enumitem}
\usepackage[utf8]{inputenc}
\usepackage{mathrsfs}
\usepackage{mathtools}
\usepackage[utf8]{inputenc} 
\usepackage{tabularx}
\usepackage{longtable}
\usepackage{makecell} 
\usepackage{booktabs} 
\usepackage[mathscr]{euscript} 
\usepackage{booktabs} 
\usepackage{multirow} 
\usepackage{amssymb,amsmath,amscd}
\usepackage{setspace}
\usepackage{color}
\usepackage[all]{xy}
\usepackage{graphics}
\usepackage[english]{babel}
\usepackage{url}
\usepackage{bm} 
\usepackage{latexsym}
\usepackage{tikz-cd}
\usepackage{fullpage}
\usepackage{thm-restate}
\usepackage{a4,color,palatino,fancyhdr}
\usepackage{graphicx}%poner [dvips] para grficos
\usepackage{float}  %Este packete permite controlar la flotacion de figuras al poner [H] o [h] como opcion de \begin{figure}
\usepackage[small, bf, margin=90pt, tableposition=bottom]{caption}
\RequirePackage{amssymb}
\RequirePackage[T1]{fontenc}
\usepackage{amscd}
\usepackage{amsrefs}
\usepackage{url}  
\usepackage[all]{xy}
\newtheorem{bigthm}{Theorem}
\newtheorem{bigcor}[bigthm]{Corollary}

\newtheorem{thm}{Theorem}[section]

\newtheorem{lem}[thm]{Lemma}

\newtheorem{prop}[thm]{Proposition}
\newtheorem{cor}[thm]{Corollary}

\theoremstyle{definition}
\newtheorem{defn}[thm]{Definition}
\newtheorem{notation}[thm]{Notation}
\newtheorem{relation}[thm]{Relation}

\theoremstyle{remark}
\newtheorem{ex}[thm]{Example}

\newtheorem{rem}[thm]{Remark}

\DeclareMathSymbol{\shortminus}{\mathbin}{AMSa}{"39}

\usepackage[colorlinks,linkcolor=blue,citecolor=blue]{hyperref}

\def\Z{\mathbb{Z}}

\def\R{\mathbb{R}}
\def\Q{\mathbb{Q}}

\def\P{\mathrm{PL}^+}
\def\T{\mathrm{Top}^{+}}
\def\coker{\mathrm{Coker}}
\def\D{\mathrm{Diff}_{\partial}}
\def\bD{\widetilde{\mathrm{Diff}}_{\partial}}

\title{On some Gromoll filtration groups and fundamental groups of $\D (D^n)$}

\author{Wei Wang}
\address{Department of Mathematics, Shanghai Ocean University, Shanghai 201306, China}
\email{weiwang@amss.ac.cn}

\date{}

\begin{document}

\maketitle

\begin{abstract}
In this paper, we compute some Gromoll filtration groups $\Gamma^{n+1}_{i+1}$ for certain $i$ when 
$8\leq n \leq 17$ and $n=4k+2\geq 18$. Using these results, we obtain some new information on $\pi_1\D (D^n)$ for $6\leq n \leq 15$ and on $\pi_2 \D (D^{4k+3})$ for $4k+3\geq 15$.
\end{abstract}
%\addtocontents{toc}{\protect\setcounter{tocdepth}{0}}
%\tableofcontents
\section{Introduction}
\subsection{The Gromoll filtration group $\Gamma_{i+1}^{n+1}$}
Let $\D(D^n)$ be the topological group of diffeomorphisms of $D^n$ which are identity near the boundary with the Whitney topology (e.g. \cite[Chapter 2]{hirsch2012differential}).

Following \cite[Section 2]{crowley2013gromoll}, let $\alpha$ be an element in $\pi_i\D(D^{n-i})$. One can represent $\alpha$ by a map
%\[
%f_{\alpha}\colon D^i\times D^{n-i}\longrightarrow D^{n-i}
%\]
%\[
%(t_1,t_2)\longmapsto f_{\alpha}(t_1,t_2),
%\]
\begin{align*}
	f_{\alpha}\colon D^i\times D^{n - i}&\longrightarrow D^{n - i}\\
	(t_1,t_2)&\longmapsto f_{\alpha}(t_1,t_2)
\end{align*}
such that $f_{\alpha}(t_1,-)\in \D(D^{n-i})$ and $f_{\alpha}$ coincides with the projection map $p_2(t_1,t_2)=t_2$ near the boundary of $D^i\times  D^{n-i}$. For $i\geq k\geq 0$, there are Gromoll homomorphisms 
\begin{align*}
\lambda^n_{k,i}\colon \pi_i\D(D^{n-i})&\longrightarrow \pi_{i-k}\D(D^{n-i+k})\\
\alpha=[f_{\alpha}]& \longmapsto \lambda^n_{k,i}(\alpha)=[f_{\lambda^n_{k,i}(\alpha)}],
\end{align*}
%\[
%\lambda^n_{k,i}\colon \pi_i\D(D^{n-i})\longrightarrow \pi_{i-k}\D(D^{n-i+k})
%\]
%\[
%\alpha=[f_{\alpha}] \longmapsto \lambda^n_{k,i}(\alpha)=[f_{\lambda^n_{k,i}(\alpha)}],
%\]
where  $\lambda^n_{k,i}(\alpha)$ is represented by a map $f_{\lambda^n_{k,i}(\alpha)}$ defined as follows
\begin{align*}
	f_{\lambda^n_{k,i}(\alpha)}\colon D^{i - k}\times D^{k}\times D^{n - i}&\cong D^{i - k}\times D^{n - i + k} \longrightarrow D^{n - i}\times D^{k}\cong D^{n - i + k}\\
(t_1,t_2,t_3)	&\longmapsto (t_1,(t_2,t_3)) \longmapsto(f_{\alpha}((t_1,t_2),t_3),t_2).
\end{align*}
%\[
%f_{\lambda^n_{k,i}(\alpha)}\colon D^{i-k}\times D^{k}\times  D^{n-i} \cong D^{i-k}\times D^{n-i+k} %\longrightarrow  D^{n-i}\times D^k\cong D^{n-i+k}
%\] 
%\[
%(t_1,t_2,t_3)\longmapsto(t_1,(t_2,t_3))\longmapsto  (f_{\alpha}((t_1,t_2),t_3), t_2).
%\]

%Let $\lambda^n_{i,i}(\pi_{i}\D (D^{n-i}))$ be the image of $\lambda^n_{i,i}$ in $\pi_0 \D(D^n)$. 
In \cite{gromoll1966differenzierbare}, Gromoll defined the group
$\Gamma_{i+1}^{n+1}$ to be the image $\lambda^n_{i,i}(\pi_{i}\D (D^{n-i}))$ in $\pi_0 \D(D^n)$
and constructed a filtration  
$$\cdots \subset \Gamma_{i+1}^{n+1}\subset \Gamma_{i}^{n+1}\subset \cdots \subset \Gamma_{1}^{n+1}.$$
Extending each element $f_{D^{n}}\in\D (D^{n})$ by the identity map to $f_{S^{n}} \in \mathrm{Diff}^+S^{n}$, the group of orientation-preserving diffeomorphisms of $S^{n}$, and clutching two $D^{n+1}$ through $f_{S^{n}}$, one can obtain a group homomorphism $ec$ from $\pi_0 \D (D^{n})$ to the group of homotopy spheres $\Theta_{n+1}$ 
\begin{align*}
	ec\colon\pi_0 \D (D^{n})&\longrightarrow \Theta_{n+1}\\
	[f_{D^{n}}]&\longmapsto [D^{n+1}\cup_{f_{S^{n}}}  D^{n+1}].
\end{align*}
When $n\geq 5$, according to \cites{cerf2006pseudo,smale1962structure}, $ec$ is an isomorphism and $\Gamma_2^{n+1}=\Gamma_{1}^{n+1}$. Furthermore, $\D(D^3)$ is contractible proved by Hatcher in \cite{hatcher1983proof}. The Gromoll filtration becomes
\[
0=\Gamma_{n-2}^{n+1}\subset \Gamma_{n-3}^{n+1}\subset \cdots \subset \Gamma_{i+1}^{n+1}\subset \cdots \subset \Gamma_{3}^{n+1}\subset \Gamma_{2}^{n+1}=\Gamma_{1}^{n+1}\cong \Theta_{n+1}.
\]
The study of this filtration is one of the central topics in geometric topology.
Many people have calculated $\Gamma_{i+1}^{n+1}$ for some values of $n$ and $i$ through various methods (see e.g. \cites{ANTONELLI1972,cerf2006pseudo,WWsurvey,crowley2013gromoll,crowley2018harmonic}). For known results on $\Gamma_{i+1}^{n+1}$, we refer to Crowley, Schick and Steimle \cites{crowley2013gromoll,crowley2018harmonic} and the references therein. 

\subsection{The relationship to $\pi_i \D(D^n)$}
The study of $\pi_i \D(D^n)$ is also a key problem in geometric topology. In recent years, fruitful breakthroughs have been made on rational homotopy groups $\pi_{i}\D(D^n)\otimes \Q$ of $\D(D^n)$ by Galatius, Krannich,
Kupers, Randal-Williams, Watanabe, Weiss and many other people.  We refer to \cite{kupers2025diffeomorphisms,randal2023diffeomorphisms,krannich2021diffeomorphisms} for more results and references. 

On the other hand, for torsion parts of $\pi_{i}\D(D^n)$, 
the Gromoll filtration groups play a role. More precisely,
let $\widetilde{\mathrm{Diff}}_{\partial}(D^n)$ be the geometric realization of the semi-simplicial group of   block diffeomorphisms of $D^n$ (see \cite[\S 1.3 \& \S 1.4]{krannich2022homological} for the definition)\footnote{According to \cite{hebestreit2021vanishing,krannich2022homological}, the definition of  semi-simplicial groups of block diffeomorphisms should satisfy a collared condition. See Remark \ref{remarkcollarcondition} for more discussions.}.
%The inclusion $\D(D^n)\hookrightarrow \widetilde{\mathrm{Diff}}_{\partial}(D^n)$ 
There is a homotopy
fibration
\[
\D(D^n)\longrightarrow \widetilde{\mathrm{Diff}}_{\partial}(D^n)\longrightarrow \frac{\widetilde{\mathrm{Diff}}_{\partial}(D^n)}{\D(D^n)},
\]
which induces long exact sequence of homotopy groups
%\[
%\cdots \longrightarrow \pi_{i}\D(D^n)\longrightarrow \pi_{i}\widetilde{\mathrm{Diff}}_{\partial}(D^n)\longrightarrow \pi_{i}\frac{\widetilde{\mathrm{Diff}}_{\partial}(D^n)}{\D(D^n)}\longrightarrow \pi_{i-1}\D(D^n)\longrightarrow \pi_{i-1}\widetilde{\mathrm{Diff}}_{\partial}(D^n)\longrightarrow \cdots
%\]
\[
\begin{tikzcd}
	\cdots \rar  &\pi_i\D(D^n) \rar &  \pi_i\widetilde{\mathrm{Diff}}_{\partial}(D^n)
	\rar  \ar[draw=none]{d}[name=X, anchor=center]{}
	&\pi_i(\frac{\widetilde{\mathrm{Diff}}_{\partial}(D^n)}{\D(D^n)} ) \ar[rounded corners,
	to path={ -- ([xshift=2ex]\tikztostart.east)
		|- (X.center) \tikztonodes
		-| ([xshift=-2ex]\tikztotarget.west)
		-- (\tikztotarget)}]{dll}[at end]{\ } \\      
	&\pi_{i-1}\D(D^n) \rar &  \pi_{i-1}\widetilde{\mathrm{Diff}}_{\partial}(D^n)
	\rar  	&\cdots.
\end{tikzcd}
\]
According to \cite[\S 2.3.2 \&  \S 2.3.3]{ANTONELLI1972} or \cite[\S 6.6]{WWsurvey}, $\pi_{i}\widetilde{\mathrm{Diff}}_{\partial}(D^n)\cong \pi_0\D(D^{n+i})\cong \Theta_{n+i+1}$ when $n+i\geq 5$, and the image of $\pi_{i}\D(D^n)$ in $\pi_{i}\widetilde{\mathrm{Diff}}_{\partial}(D^n)\cong \pi_0\D(D^{n+i})$ coincides with the Gromoll filtration group
$\Gamma^{n+i+1}_{i+1}$, i.e. the following diagram is commutative
\[
\begin{CD}
\pi_{i}\D(D^n) @>>>  \pi_i \bD (D^n)\\
@V \lambda^{n+i}_{i,i}VV     @V \cong VV \\
\pi_0 \D(D^{n+i})   @> ec >>   \Theta_{n+i+1}.  
\end{CD}
\]
It follows that one has the short exact sequence

\begin{equation}\label{maineq}
0\longrightarrow \Gamma^{n+i+2}_2/\Gamma^{n+i+2}_{i+2}\longrightarrow \pi_{i+1}\frac{\widetilde{\mathrm{Diff}}_{\partial}(D^n)}{\D(D^n)}\longrightarrow \pi_{i}\D(D^n)\longrightarrow
\Gamma^{n+i+1}_{i+1}\longrightarrow 0,
\end{equation}
where the Gromoll filtration groups take affects if we consider the torsion part of $\pi_{i}\D(D^n)$.

\subsection{Main results}
%Our main results are
%\begin{thm} Let $\Gamma_{i+1}^{n+1}$ be the Gromoll filtration group.
%	\begin{itemize}
%		\item[(1)] $\Z_3\subset \Gamma^{10}_{5}\subset \Gamma_{2}^{10}\cong \Z_2\oplus \Z_3$\\
%		\item[(2)] $\Z_8\subset \Gamma^{18}_{10}\subset \Gamma_2^{18}\cong \Z_2 \oplus \Z_8$
%	\end{itemize}
%\end{thm}
%According to \cite{Crowley2}, $\Z_2\oplus 0\subset \Gamma^{10}_{4}$ and $\Z_2\oplus 0\subset \Gamma_{14}^{18}$. We have
%\begin{coro}
%	\begin{itemize}
%		\item[(1)]  $\Gamma^{10}_{i+1}=\Gamma^{10}_{2}\cong \Z_2\oplus \Z_3$ for $i\leq 4$.\\
%		 \item[(2)]  $\Gamma^{18}_{i+1}=\Gamma^{18}_{2}\cong \Z_2\oplus \Z_3$ for $i\leq 14$.\\
%	\end{itemize}
%\end{coro}

\subsubsection{Gromoll filtration groups}
Our first main result is the calculation of some Gromoll filtration groups.
\begin{bigthm}\label{Theorem A}
	Let $\Gamma_{i+1}^{n+1}$ be the $(i+1)$-th Gromoll filtration group of $\Theta_{n+1}$
	\begin{itemize}
		\item[(a)] $\Gamma^{10}_{4}=\Gamma_{2}^{10}$, $\Gamma_{7}^{13}=\Gamma_{2}^{13}$, and $\Gamma_{5}^{18}=\Gamma_{2}^{18}$.
		\item[(b)] $(\Z_{2})^2\subset \Gamma_{3}^9$, $\Z_{2}\oplus 0\subset \Gamma_9^{15}$, and $(\Z_{2})^3\subset \Gamma_{5}^{17}$.
		\item[(c)] $\Gamma_{4}^{4k+3}=\Gamma_{3}^{4k+3}$ with $\Theta_{4k+3}/\Gamma_{4}^{4k+3}\cong \Z_{2}$ when $k\geq 3$.
	\end{itemize}
\end{bigthm}

As a supplement to the tables in \cite[Appendix A]{crowley2013gromoll} and \cite[Appendix A]{crowley2018harmonic}, we update certain Gromoll filtration groups in Table \ref{tableGromoll} of Appendix \ref{appdixA} by assembling results from Theorem \ref{Theorem A} and the literature of the Gromoll filtration groups.

\subsubsection{Homotopy groups of $\D(D^n)$}
The information of the Gromoll filtration groups in Theorem \ref{Theorem A} helps us to obtain some information on $\pi_i \D(D^n)$. We first discuss the fundamental group of $\D(D^n)$ when $6\leq n \leq 15$. 
\begin{bigthm}\label{thmpi1}
	\begin{itemize}
		\item[(a)] The homomorphism 
	$$\pi_{1}\D(D^n)\rightarrow \pi_1 \bD(D^n)\cong \Theta_{n+2}$$ 
	is split surjective when $6\leq n \leq 15$ and $n\neq 9,13$. \\
	\item[(b)] The composition 	
	$$\pi_{1}\D(D^{13})\rightarrow \pi_1 \bD(D^{13})\cong \Theta_{15}\rightarrow \Theta_{15}/bP_{16}\cong \Z_2$$ 
	is split surjective.
	\end{itemize}
\end{bigthm}
Based on Theorem \ref{thmpi1} and the exact sequence (\ref{maineq}), we can show
\begin{bigcor}\label{CorollaryC}
\begin{itemize}
%	\item[(a)] $\pi_{1} \D(D^{10})\cong \Z_{2}$, $\pi_{1}\D(D^{12})\cong \Z_{2}$, and $\pi_{1}\D (D^{15})\cong \Z_{2}\oplus \Theta_{17}$.
	\item[(a)] $\pi_{1}\D (D^{15})\cong \Z_{2}\oplus \Theta_{17}$.\\
		\item[(b)] $ \pi_{1} \D (D^7)\cong \Theta_{9}\oplus \pi_2(\frac{\widetilde{\mathrm{Diff}}_{\partial}(D^7)}{\D(D^7)} )$, and $\pi_{1}\D(D^{9})/\Z_{2} \cong \Theta_{11}$.\\
	
%	\item[(c)] When $n=6,8,11,14$, there exists subgroups in $\pi_{1}\D(D^n)$ maps one-to-one onto $\Theta_{n+2}$ through the homomorphism $\pi_{1}\D(D^n)\rightarrow \pi_1 \bD(D^n)\cong \Theta_{n+2}$
%	\[
%	\Z_2\subset \pi_{1}\D(D^6), \Z_6	\subset \pi_{1}\D(D^8), 	\Z_3\subset \pi_{1}\D(D^{11}), \  \text{and}\ 	%\Z_2\subset \pi_{1}\D(D^{14}).
%	\]
\end{itemize}
\end{bigcor}

Some information on the groups $\pi_1\D (D^n)$ and some of their generators is collected in Table \ref{tableDiff} of Appendix \ref{appdixB}. 
We also have some results on $\pi_{2}\D(D^{4k-1})$ when $k\geq 4$.

%\begin{bigthm}\label{Theorem B}
%	The following homomorphisms are split surjective
%	\begin{itemize}
%			\item[(a)] $\pi_{i}\D(D^6)\longrightarrow \pi_i \bD(D^6)$ when $i=1,3,6;$\\
%		$\pi_{i}\D(D^6)\rightarrow \pi_i \bD(D^6)\rightarrow \Theta_{6+i+1}/bP_{6+i+2}$ when $i=2,8$.
%			\item[(b)] $\pi_{i}\D(D^{14})\longrightarrow \pi_i \bD(D^{14})$
%			           when $i=1,3$.
%			\end{itemize}
%\end{bigthm}
\begin{bigcor}\label{CorollaryD}
When $k\geq 4$, there is an isomorphism
$$\pi_{2}\D(D^{4k-1})\cong \Gamma_3^{4k+2}.$$
Moreover, $\pi_{2}\D(D^{15})\cong \Gamma_3^{18}\cong \Theta_{18}$.
\end{bigcor}
More generally, we can show
\begin{bigcor}\label{CorollaryE}
	Let $M^{4k-1}$ be a 3-connected compact smooth manifold with $k\geq 4$, the homomorphism 
	\[
	\pi_{2}\D(M^{4k-1})\longrightarrow \pi_2 \bD(M^{4k-1})
	\]
	is injective with order $|\pi_2 \bD(M^{4k-1})/\pi_{2}\D(M^{4k-1})|\leq 2$.
\end{bigcor}

%\subsubsection{Fundamental groups of diffeomorphism groups of homotopy spheres}
%\subsubsection{Realization problem of $\mathrm{Diff}S^n$}
\subsection{Strategy and organization}
Let $\alpha\in \pi_i SO(j)$ and $\beta\in \pi_jSO(i)$, Antonelli, Burghelea and Kahn \cite[Page 15]{ANTONELLI1972} proved that the Milnor--Novikov pairing $\sigma_{i,j}(\alpha\otimes\beta)$ of $\alpha$ and $\beta$
\begin{align*}
	\sigma_{i,j}\colon \pi_i SO(j)\otimes \pi_j SO(i)&\longrightarrow \pi_0 \D(D^{i+j})\cong \Theta_{i+j+1} \\
	\alpha\otimes\beta&\longmapsto \sigma_{i,j}(\alpha\otimes\beta)
\end{align*}
can be lifted to $\pi_{a+b+1}\D(D^{i+j-a-b-1})$ and the element $\sigma_{i,j}(\alpha\otimes\beta)$ lies in the Gromoll filtration group $\Gamma_{a+b+2}^{i+j+1}$ when $\alpha$ and $\beta$ come from $\pi_i SO(j-a)$ and $\pi_jSO(i-b)$ through the natural inclusions.
On the other hand, Frank \cite{frank1968invariant} and Smith \cite{smith1974framings} found that the quotient class of  $[\sigma_{i,j}(\alpha\otimes\beta)]\in \Theta_{i+j+1}/bP_{i+j+2}\subset \mathrm{Coker}J_{i+j+1}$  belongs to a quotient set of a certain Toda bracket when $\alpha$ and $\beta$ satisfy some Toda type conditions. With this inclusion, the element $\sigma_{i,j}(\alpha\otimes\beta)$ can be determined in some dimensions (see Section \ref{MNparingToda} for more details). However, the proof stated in \cite{frank1968invariant} lacked several details. Following \cite{FarrellHsiang1976} and \cite{frank1967reducible}, we provide a more detailed proof in Section \ref{SectionProofofFrankSmith}.

Our strategy is to combine these two approaches. We use the results of Antonelli, Burghelea and Kahn (Theorem \ref{theoremABK}) to construct the element $\sigma_{i,j}(\alpha\otimes\beta)$ in $\Gamma_{a+b+1}^{i+j+1}$, and apply the Frank--Kosinski--Smith theorem (Theorem \ref{FrankSmith}) to determine its image in $\mathrm{Coker}J_{i+j+1}$. We also find that, in certain cases, the lifting of $\sigma_{i,j}(\alpha\otimes\beta)$ enables the subgroup generated by $\sigma_{i,j}(\alpha\otimes\beta)$ to split via the Gromoll homomorphism. 

Once we have obtained information about certain Gromoll filtration groups via the aforementioned approach, we can use the short exact sequence (\ref{maineq}) to get further information about $\pi_1 \D(D^n)$. This relies on the Hatcher spectral sequence and the commutative composition diagrams used by Burghelea and Lashof \cite{BurgheleaLashof}, as well as by Crowley, Schick and Steimle  \cite{crowley2013gromoll,crowley2018harmonic}.  

This paper is organized as follows. In Section \ref{SectionMNTodaCompoistion}, we recall two definitions of the Milnor--Novikov pairing, as well as Antonelli, Burghelea and Kahn's theorem, the Frank--Kosinski--Smith theorem, and the commutative composition diagram. In Section \ref{SectionMilnor}, we construct some typical Milnor--Novikov pairings using the methods introduced in Section \ref{SectionMNTodaCompoistion}.
We prove Theorem \ref{Theorem A}
and Theorem \ref{thmpi1} in Section \ref{SectionProofMainthm}, and we also give the proofs of Corollaries \ref{CorollaryC}, \ref{CorollaryD}, and \ref{CorollaryE} as applications. 
A detailed proof of the Frank--Kosinski--Smith theorem is provided in Section \ref{SectionProofofFrankSmith}.
%\subsection*{Acknowledgments}
%We would like to thank Manuel Krannich and Cary Malkiewich for their helpful discussions.
\subsection{Acknowledgments}
The author would like to thank the referee for their numerous helpful comments and
suggestions, which have improved this paper in many aspects. In particular, the author owes much gratitude to the referee for refining the statement of Theorem \ref{FrankSmith} and improving Remark \ref{remarkofFKSthm}. The author benefited greatly from discussions with Lizhen Qin, Manuel Krannich, 
Cary Malkiewich, Yang Su, and Xiaolei Wu. The author would also like to thank Jie Zhang for assistance with the references.

\section{Milnor--Novikov pairing, Toda bracket and composition commutative diagram}\label{SectionMNTodaCompoistion}
The Milnor--Novikov pairing is both a classical method to construct homotopy spheres and a method to construct non-trivial elements in $\pi_i \D(D^n)$. In this section, we introduce the definition of the Milnor--Novikov pairing from two points of view. From these viewpoints, Antonelli, Burghelea and Kahn's theorem (Theorem \ref{theoremABK}) relates this paring to the Gromoll filtration groups.% and the homotopy groups of$\D(D^n)$.  

On the other hand, the Frank--Kosinski--Smith theorem (Theorem \ref{FrankSmith}) provides a way to connect certain Milnor--Novikov parings to Toda brackets, some 
of which yield homotopy spheres that map non-trivially to the cokernel of the $J$-homomorphism.

\subsection{The Milnor--Novikov pairing}
\subsubsection{Plumbing construction} Let $\mathrm{Vect}_{S^{k+1}}(\R^{m})$ be the set of isomorphism classes of oriented $m$-dimensional euclidean vector bundles over $S^{k+1}$. There is a one-to-one correspondence between the set $\mathrm{Vect}_{S^{k+1}}(\R^m)$ and the homotopy group $\pi_{k}SO(m)$
\begin{align*}
	\pi_{k}SO(m)& \longrightarrow \mathrm{Vect}_{S^{k+1}}(\R^m)\\
	\gamma &\longmapsto \xi_{\gamma}.
\end{align*}

For $\xi_1\in \mathrm{Vect}_{S^{i+1}}(\R^{j+1})$ and $\xi_2\in \mathrm{Vect}_{S^{j+1}}(\R^{i+1})$,
denote by $D(\xi_{1})$ and $D(\xi_{2})$ the disc bundles of $\xi_{1}$ and $\xi_{2}$. Let $D(\xi_{1})\cup_{\mathrm{Plumb}}D(\xi_{2})$ be the plumbing of $D(\xi_{1})$ and $D(\xi_{2})$ by gluing 
two coordinate chars $D^{i+1}\times D^{j+1}\subset D(\xi_{1})$ and $D^{j+1}\times D^{i+1}\subset D(\xi_{2})$ to each other via interchanging their factors. $D(\xi_{1})\cup_{\mathrm{Plumb}}D(\xi_{2})$ can be a compact smooth manifold by smoothing corners
(see e.g. \cite{browder1966open,crowley2017positive,milnor1959differentiable} 
for a more detailed definition of the plumbing construction). 

Let $s_1\colon \pi_* SO(n)\longrightarrow \pi_* SO(n+1)$
be the homomorphism induced by the standard inclusion $SO(n)\hookrightarrow SO(n+1)$.
For elements $\alpha\in \pi_iSO(j)$ and $\beta\in \pi_jSO(i)$, let $\xi_{s_1 \alpha}\in \mathrm{Vect}_{S^{i+1}}(\R^{j+1})$ and $\xi_{s_1\beta}\in \mathrm{Vect}_{S^{j+1}}(\R^{i+1})$ be two euclidean vector bundles corresponding to $s_1\alpha\in \pi_i SO(j+1)$ and $s_1\beta\in \pi_j SO(i+1)$. It is known that (\cite[Lemma 1, Page 963]{milnor1959differentiable}) the boundary $\partial (D(\xi_{s_1\alpha})\cup_{\mathrm{Plumb}}D(\xi_{s_1\beta}))$ of the plumbing $D(\xi_{s_1\alpha})\cup_{\mathrm{Plumb}}D(\xi_{s_1\beta})$ is a smooth homotopy sphere. According to \cite[Page 834]{Lawson1973remarks}, we have the following bilinear pairing
\begin{align*}
\sigma_{i,j}\colon  \pi_iSO(j)\otimes \pi_j SO(i) &\longrightarrow \Theta_{i+j+1}\\
                   (\alpha,\beta)                 &\longmapsto [\partial (D(\xi_{s_1\alpha})\cup_{\mathrm{Plumb}}D(\xi_{s_1\beta}))].\\
\end{align*}

\subsubsection{Commutator construction}
There is an another construction of the pairing $\sigma_{i,j}$ which relates closely to Gromoll filtration groups (\cite[Section 1.2]{ANTONELLI1972}).
For $\alpha\in \pi_iSO(j)$ and $\beta\in \pi_jSO(i)$, we can represent $\alpha$ and $\beta$ by
smooth maps 
\[
f_{\alpha}\colon D^i\longrightarrow SO(j),f_{\beta}\colon D^j\longrightarrow SO(i),
\]
such that $f_{\alpha}(x)$ and $f_{\beta}(y)$ are identity matrices when $x$ and $y$ are in certain small collars of $D^i$ and $D^j$.

Consider two diffeomorphisms $A_{\alpha},B_{\beta}\colon D^i\times  D^j\longrightarrow D^i\times D^j$ induced by 
$f_{\alpha}$ and $f_{\beta}$ with
\begin{align*}
A_{\alpha}\colon (x,y)\longmapsto (x,f_{\alpha}(x)y), \\
 B_{\beta}\colon (x,y)\longmapsto (f_{\beta}(y)x,y),
\end{align*}
the commutator $A^{-1}_{\alpha}B^{-1}_{\beta}A_{\alpha}B_{\beta}$ is a diffeomorphism of $D^i\times D^j$ which is identity near the boundary of $D^i\times D^j$, i.e.  $A^{-1}_{\alpha}B^{-1}_{\beta}A_{\alpha}B_{\beta}\in \D (D^{i+j})$. Therefore, we have the following pairing
\begin{align*}
\sigma^c_{i,j}\colon\pi_iSO(j)\otimes \pi_j SO(i)&\longrightarrow \pi_0 \D (D^{i+j})\\
	(\alpha,\beta)&\longmapsto [A^{-1}_{\alpha}B^{-1}_{\beta}A_{\alpha}B_{\beta}].
\end{align*}
Recall that when $i+j\geq 5$, $ec\colon\pi_0 \D (D^{i+j})\longrightarrow \Theta_{i+j+1}$ is an isomorphism and one can identify $\sigma_{i,j}$ with $\sigma^c_{i,j}$ through the following commutative diagram (\cite[Page 834]{Lawson1973remarks})
\[
\begin{CD}
		\sigma^c_{i,j}\colon @. \pi_iSO(j)\otimes \pi_j SO(i) @>>> \pi_0 \D (D^{i+j}) \\
		   @.         @|                                               @VV ec V\\
		  \sigma_{i,j}\colon @. \pi_iSO(j)\otimes \pi_j SO(i) @>>>    \Theta_{i+j+1}.      
\end{CD}
\]

\subsection{Milnor--Novikov pairings and Gromoll filtration groups}
In \cite[Page 15]{ANTONELLI1972},  Antonelli, Burghelea and Kahn found two important relations between Gromoll filtration groups and Milnor--Novikov pairings through  the above commutator constructions. They proved (\cite[Proposition 1.2.4]{ANTONELLI1972})
\begin{thm}\label{theoremABK}
For any integers $a,b$ with $0\leq a \leq j$ and $0\leq b \leq i$, let 
\[
s_a:\pi_{i} SO(j-a)\longrightarrow \pi_i SO(j) \text{\quad and \quad} s_b:\pi_{j} SO(i-b)\longrightarrow \pi_{j}SO(i)
\]
be the homomorphisms induced by the standard inclusions $SO(j-a)\longrightarrow SO(j)$ and  $SO(i-b)\longrightarrow SO(i).$ 
\begin{itemize}
	\item[(a)] Denote $c=a+b+1$, there is a homomorphism
 $$\sigma_{i,j}^{(a,b)}\colon \pi_iSO(j-a)\otimes \pi_j SO(i-b)\longrightarrow \pi_{c} \D (D^{i+j-c})$$
which makes the following diagram commutative
\[
\begin{CD}
	 \pi_iSO(j-a) @.\otimes @.\pi_j SO(i-b) @>\sigma_{i,j}^{(a,b)}>>  \pi_{c}\D (D^{i+j-c})\\
         @Vs_aVV         @.           @Vs_bVV                @VV \lambda_{c,c}^{i+j} V\\
	 \pi_iSO(j)   @.\otimes @. \pi_j SO(i) @>\sigma_{i,j}>>     \pi_0 \D ( D^{i+j}).      
\end{CD}
\]
\item[(b)]$\sigma_{i,j}(s_a (\pi_iSO(j-a)) \otimes s_b(\pi_j SO(i-b) ) )\subseteq \Gamma^{i+j+1}_{c+1}$.
\end{itemize}
\end{thm}

\subsection{Milnor--Novikov parings and Toda brackets}\label{MNparingToda}
Denote by $\pi^s_i$ the $i$-th stable homotopy group of spheres.
Write $\Sigma\colon\pi_{n+i} S^n\longrightarrow \pi_{n+i+1}S^{n+1}$ for the suspension homomorphism and  $\Sigma^{\infty}\colon$ $\pi_{n+i} S^n\longrightarrow \pi_{i}^{s}$ for the stabilisation.
Suppose $i>0$, let $J_{i,n}\colon \pi_i SO(n)\cong \pi_iO(n)\longrightarrow \pi_{i+n} S^n$ be the unstable $J$-homomorphism and 
$J_{i}\colon \pi_i SO\cong \pi_i O\longrightarrow \pi_{i}^s$ be the stable $J$-homomorphism.

Denote by $bP_{i+1}$ the subgroup of $\Theta_i$ consisting of homotopy $i$-spheres which bounds parallelizable manifolds.
One has the homomorphism $\Theta_{i}/bP_{i+1}\longrightarrow \mathrm{Coker} J_i$ which is an isomorphism when $i\not\equiv 2 ( \mathrm{mod}\ 4)$ and is injective with index at most 2 when $i\equiv 2 (\mathrm{mod}\ 4)$ (\cite[Theorem 12.1]{luck2024surgery}).

Let $P$ be the composition of $\Theta_i\longrightarrow \Theta_{i}/bP_{i+1}\longrightarrow \mathrm{Coker} J_i$ and $q$ be the quotient homomorphism $\pi_i^s \longrightarrow \coker J_i$. 
The following theorem relates Toda brackets to Milnor--Novikov pairings.

\begin{restatable}{thm}{FrankKosinskiSmiththm}\textup{(The Frank--Kosinski--Smith Theorem).}\label{FrankSmith}
For $\alpha\in \pi_i SO(j)$ and $\beta \in \pi_j SO(i)$ with $i,j\geq 1$, let $u=J_{j,i}\beta$ and suppose
\begin{itemize}

	\item[(1)] $J_{i,j}\alpha= vw$ with $v$ finite order in $\pi_{j+t}S^j$,  $ 0<t<i$;
	\item[(2)] 	$s_{\infty}\alpha=0$, where $s_{\infty}\colon \pi_i SO(j)\longrightarrow \pi_i SO$ is the homomorphism induced by the natural inclusion;
	% $(\Sigma^{\infty}v) (\Sigma^{\infty}w)=0$;

	\item[(3)] $(\Sigma^{\infty}u) (\Sigma^{\infty}v)=0$,
%	\item[(4)] $(\Sigma^{\infty}u)\circ \pi^s_{i+1}\subset \mathrm{Im}J_{i+j+1}$,
\end{itemize}
then one has the following inclusion
\[
	(-1)^{(i+1)(j+1)-1}P(\sigma_{i,j}(\alpha,\beta))\in q(<\Sigma^{\infty} u, \Sigma^{\infty}v, \Sigma^{\infty}w>),
\]
where $<\Sigma^{\infty}u, \Sigma^{\infty}v, \Sigma^{\infty}w>$ denotes the stable Toda bracket of $\Sigma^{\infty} u$, $\Sigma^{\infty}v$ and $\Sigma^{\infty}w$.
\end{restatable}

\begin{rem}\label{remarkofFKSthm}
The statement of Theorem \ref{FrankSmith} differs from that of Theorem 2 in \cite{frank1968invariant}, which involves the notation "twisted suspension" which Frank did not define in \cite{frank1968invariant}. We believe that Frank's version of Theorem 2 in \cite{frank1968invariant} is more general, while in Theorem \ref{FrankSmith} we have imposed certain technical conditions to avoid introducing  the notion of "twisted suspension" for simplicity.

When $J_{i,j}(\alpha)$ decomposes into elements in the image of the $J$-homomorphism, Kosinski gave a full proof under the dimension assumptions he makes (\cite[Theorem 4.3 (b)]{kosinski1971toda}). In \cite{smith1974framings}, Smith stated a version in Theorem 5.2 without dimension assumptions, but his proof only refers to that in \cite{kosinski1971toda} and does not provide details on how these dimension restrictions are removed.

We provide details of the proof of Theorem \ref{FrankSmith} in Section \ref{SectionProofofFrankSmith}, following the ideas in \cite{frank1967reducible,frank1968invariant}.

\end{rem}
\subsection{Morlet equivalences and composition diagrams}
Let $\mathrm{PL}_n$ be the simplicial group of piecewise linear homeomorphisms of $\R^n$ fixing the origin, 
and let $B\mathrm{PL}_n$ be its classifying space which classifies piecewise linear $\R^n$-bundles. Denote $\mathrm{PL}=\lim\limits_{n\to\infty} \mathrm{PL}_n $ and $B\mathrm{PL}=\lim\limits_{n \to \infty} B\mathrm{PL}_n $. 

We also write $\mathrm{PL}_n/O_n$ for the homotopy fiber of the natural map $BO(n)\longrightarrow B \mathrm{PL}_n$ and write $\mathrm{PL}/O$ for the homotopy fiber of  $BO\longrightarrow B \mathrm{PL}$. The Morlet equivalence $\D(D^n)\simeq \Omega^{n+1} \mathrm{PL}_n/O_n$ (\cite[Theorem 4.4 (a)]{BurgheleaLashof}) induces an isomorphism
\[
M_n\colon \pi_i \D(D^n)\longrightarrow \pi_{n+i+1}\mathrm{PL}_n/O_n.
\]

According to \cite[$\S$ 7]{BurgheleaLashof}, for any pair $(k,t)$ with $t\leq k-2$, precomposition defines a group homomorphism
$\tilde{\psi}^X_{k,t}\colon \pi_k(X)\otimes \pi_t^s\longrightarrow \pi_{k+t}(X)$
for any pointed topological space $X$. We define $$\psi^n_{i;t}\colon 	\pi_i \D (D^n) \otimes \pi^s_t\longrightarrow \pi_{i+t} \D (D^n) $$ to be the composition 
\[
\begin{CD}
	\pi_i \D (D^n)  \otimes \pi^s_t @>M_n\otimes id >>\pi_{n+i+1} \mathrm{PL}_n/O_n \otimes \pi^s_t @>\tilde{\psi}^{\mathrm{PL}_n/O_n}_{n+i+1,t} >>  \pi_{n+i+1+t} \mathrm{PL}_n/O_n@> M_n^{-1}>>	\pi_{i+t} \D (D^n).
\end{CD}
\]
By \cite[Theorem 7.1]{BurgheleaLashof}, when $t\leq n+i-1$, the following commutative diagram holds
\[
\begin{CD}
@.	\pi_i \D (D^n) \otimes \pi^s_t @> \psi^n_{i;t} >> 	\pi_{i+t} \D (D^n) \\
@.	@V M_n VV                             @V M_n VV\\
@.	\pi_{n+i+1} \mathrm{PL}_n/O_n \otimes \pi^s_t  @> \tilde{\psi}^{\mathrm{PL}_n/O_n}_{n+i+1,t} >> 	\pi_{n+i+1+t} \mathrm{PL}_n/O_n\\
@.	@V S_* VV                                 @V S_* VV\\
\Theta_{n+i+1}\otimes \pi^s_t @> \cong >>		\pi_{n+i+1} \mathrm{PL}/O \otimes \pi^s_t  @> \tilde{\psi}^{\mathrm{PL}/O}_{n+i+1,t} >> 	\pi_{n+i+1+t} \mathrm{PL}/O\\
@VVV			@VVV                                 @VVV\\
\mathrm{Coker} J_{n+i+1}\otimes \pi^s_t @>>>	\pi_{n+i+1} \mathrm{G}/O \otimes \pi^s_t  @> \tilde{\psi}^{\mathrm{G}/O}_{n+i+1,t} >> 	\pi_{n+i+1+t} \mathrm{G}/O,\\
		\end{CD}
\]
where $S_*$ denotes the homomorphism induced by the natural inclusion, and $G=\lim\limits_{n \to \infty}G_n$ denotes the homotopy direct limit of the spaces $G_n$ of homotopy equivalences of $S^n$. 

Meanwhile, according to \cite[Lemma 2.5]{crowley2013gromoll} and \cite[Page 1082]{crowley2018harmonic}, the Gromoll homomorphism $\lambda^{n+i}_{i,i}$ and the Morlet equivalence $M_n$ commute as follows
\[
\begin{CD}
	\pi_i \D (D^n) @> M_n >>    	\pi_{n+i+1} \mathrm{PL}_n/O_n\\
	@V \lambda^{n+i}_{i,i} VV    @VS_*VV\\
	\pi_0 \D(D^{n+i}) @>M_{n+i}>> 	\pi_{n+i+1} \mathrm{PL}_{n+i}/O_{n+i}.
\end{CD}
\]
Combining these results, we obtain
\begin{prop}\label{Compsitiondiagram}
When $t\leq n+i-1$,	the following commutative diagram holds
		$$	\begin{CD}
			\pi_i \D (D^n) \otimes \pi^s_t @> \psi^n_{i;t}>> 	\pi_{i+t} \D (D^n) \\
			@V \lambda^{n+i}_{i,i}\otimes id VV                                @V\lambda^{n+i+t}_{i+t,i+t} VV\\
			\pi_0 \D(D^{n+i})\otimes \pi^s_t  @> \hat{\psi}^{n+i}_{0;t}>>   \pi_0 \D (D^{n+i+t})\\
			@V P\otimes id VV       @V P VV\\
			\mathrm{Coker} J_{n+i+1}\otimes \pi^s_t @>precomposition >> 	\mathrm{Coker} J_{n+i+t+1},
		\end{CD}$$
		where $\hat{\psi}^{n+i}_{0;t}$ denotes the composition of $\lambda^{n+i+t}_{t,t}$ and $\psi^{n+i}_{0;t}$.
		\end{prop}
		\begin{proof}
		We only mention that $\mathrm{Im} J_{n+i+1}\circ\pi^s_t\subset \mathrm{Im} J_{n+i+1+t}$ when $t\leq n+i-1$. 
		\end{proof}
		\begin{cor}\label{Gromollcomposition}
	For the Gromoll filtration groups, one has the following commutative diagram
	\[
	\begin{CD}
		\Gamma_{i+1}^{n+i+1} \otimes \pi^s_t  @> \hat{\psi}^{n+i}_{0;t} >>  	\Gamma_{i+t+1}^{n+i+t+1}\\
		@V P\otimes id VV       @V P VV\\
	\mathrm{Coker} J_{n+i+1}\otimes \pi^s_t @>Precomposition >> 	\mathrm{Coker}J_{n+i+t+1}.
	\end{CD}
	\]
		\end{cor}
\begin{proof}
The Gromoll filtration groups $\Gamma_{i+1}^{n+i+1}$ and $\Gamma_{i+t+1}^{n+i+t+1}$ are images of $ \lambda^{n+i}_{i,i}$ and $\lambda^{n+i+t}_{i+t,i+t}$. We get this commutative diagram of Gromoll filtration groups by restricting the diagram of Proposition \ref{Compsitiondiagram} to the images of  $ \lambda^{n+i}_{i,i}$ and $\lambda^{n+i+t}_{i+t,i+t}$.

\end{proof}

\section{Some Milnor--Novikov pairings}\label{SectionMilnor}

In this section, we use Theorem \ref{FrankSmith} to construct some Milnor--Novikov pairings whose images are nontrivial in $\mathrm{Coker}J$. We also construct some nontrivial elements in $\pi_i \D(D^{6})$ and $\pi_i \D(D^{14})$ for certain $i$ via Theorem \ref{theoremABK} and Proposition \ref{Compsitiondiagram}.

In this section, lots of symbols and data on the homotopy groups of $S^n$, $SO(n)$ and the groups of homotopy sphere $\Theta_n$ are used. For convenience, we list them in Appendix \ref{Appendixhomotopyinformation}.
% based on some examples in \cite{bier2006detecting,frank1968invariant,stolz2007hochzusammenhangende}.
\subsection{Dimension 8 and 16} 
We first consider two pairings in dimension 8 and 16, which were also discussed by Bier and Ray in {\cite{bier2006detecting}}, by
Frank in \cite[Example 1, Page 565]{frank1968invariant}, and by Stolz in {\cite[Satz 12.1]{stolz2007hochzusammenhangende}}. 

In dimension 8, by Proposition \ref{unstableJ}, the unstable $J$-homomorphism 
$$J_{3,3}\colon \pi_3 SO(3)\longrightarrow \pi_6 S^3$$
is surjective. Choose an element $\tilde{t}_3\in \pi_3SO(3)$ with $J_{3,3}(\tilde{t}_3)=\nu'\in \pi_6S^3$. Let $\eta_3\in \pi_4 S^3$ be the generator and $\gamma_{3,4}\in \pi_3 SO(4)$ be the second generator with $J_3 s_{\infty}\gamma_{3,4}=J_3\gamma_3=\nu +\alpha_1\in\pi^s_3$. 

For the following Milnor--Novikov pairing
\[
\begin{CD}
	\pi_4SO(3) @.\otimes @.\pi_3 SO(4) @>\sigma_{4,3}^{(0,0)}>>  \pi_{1}\D (D^{6})\\
	@Vs_0VV         @.           @Vs_0VV                @VV \lambda_{1,1}^{7} V\\
	\pi_4SO(3)   @.\otimes @. \pi_3 SO(4) @>\sigma_{4,3}>>     \pi_0 \D (D^{7})\cong \Theta_8,     
\end{CD}
\]
let $\theta_1^6=\sigma_{4,3}^{(0,0)}(\tilde{t}_3\eta_3,9\gamma_{3,4})\in \pi_{1}\D (D^{6})$, one has

\begin{prop}\label{dim8z2}

\begin{itemize}
		\item[(a)] $\lambda_{1,1}^{7}(\theta_1^6)=\Sigma^8_{[\epsilon]}\in \Theta_8$;
		\item[(b)] $\theta_1^6$ is of order 2;
		\item[(c)] the forgetful homomorphism $ \pi_{1} \D (D^6)  \longrightarrow \pi_1 \bD(D^6) $ is split surjective.
\end{itemize}
\end{prop}
\begin{proof}
For convenience, we choose $9\gamma_{3,4}$ to cancel the 3-torsion. Write $\alpha$ for $\tilde{t}_3\eta_3$ and $\beta$ for $9\gamma_{3,4}$, note that
\begin{itemize}
	\item[(1)]  $J_{4,3}(\alpha)=J_{4,3}(\tilde{t}_3\eta_3)=J_{3,3}(\tilde{t}_3 )\eta_6=\nu'\eta_6$ with  $\Sigma^{\infty}\nu'=2\nu$ and $\Sigma^{\infty}\eta_6=\eta$;
		\item[(2)] $\pi_4 SO=0$ implies $s_{\infty}\alpha=0$;
	\item[(3)] $\Sigma^{\infty}J_{3,4}(\beta)=9(\nu+\alpha_{1})=\nu$;
	\item[(4)]  The elements $\nu,2\nu$ and $\eta$ satisfy the Toda condition.
\end{itemize}
By Theorem \ref{FrankSmith}, we have
\[
P\sigma_{4,3}(\alpha,\beta)=P\sigma_{4,3}(\tilde{t}_3\eta_3,9\gamma_{3,4})\in -q(<\nu,2\nu,\eta>)=-[\epsilon]=[\epsilon].
\]
In dimension 8, the composition map $P\colon \Theta_8\cong\Z_2 \Sigma^8_{[\epsilon]}\longrightarrow \mathrm{Coker}J_8\cong \Z_2[\epsilon]$ is an isomorphism and we have 
$$\lambda_{1,1}^{7}(\theta_1^6)=\lambda_{1,1}^{7}(\sigma_{4,3}^{(0,0)}(\tilde{t}_3\eta_3,9\gamma_{3,4})))=\sigma_{4,3}(\tilde{t}_3\eta_3,9\gamma_{3,4}))=\Sigma^8_{[\epsilon]}.$$ 
The order of $\tilde{t}_3\eta_3$ is 2,  which implies that the element $\theta_1^6$ is of order 2. Moreover,
the subgroup $\Z_{2} \theta_1^6\subset \pi_{1} \D (D^6)$ fits into the following commutative diagram
\[
\begin{CD}
	\Z_{2} \theta_1^6@> \cong >> \Theta_{8}\\
	@VVV       @A \cong AA\\
	\pi_{1} \D (D^6)  @>>> \pi_1 \bD(D^6).
\end{CD}
\] 
Then the homomorphism $ \pi_{1} \D (D^6)  \longrightarrow \pi_1 \bD(D^6) $ is split surjective.

\end{proof}
\begin{cor}\label{corodim9}
The composition	$ \pi_{2} \D (D^6)  \longrightarrow \pi_2 \bD(D^6)\cong \Theta_{9}\longrightarrow \Theta_9/bP_{10}$ is split surjective.
\end{cor}
\begin{proof}
 By \cite[Theorem 1.1]{crowley2018harmonic}, choose the $2$-torsion element $\tau_2^6\in \pi_2 \D (D^6)$ with $\alpha$-invariant 1 (for the definition of the $\alpha$-invariant, see e.g. \cite[Section 2.3]{crowley2013gromoll}). Let $\eta\in \pi^s_1$ be the generator. By Proposition \ref{Compsitiondiagram}, the following composition diagram
 $$	\begin{CD}
 	\pi_1 \D (D^6) \otimes \pi^s_1 @> \psi^6_{1;1}>> 	\pi_{2} \D (D^6) \\
 	@V \lambda^7_{1,1}\otimes id VV                                @V\lambda^8_{2,2} VV\\
 	\pi_0 \D(D^{7})\otimes \pi^s_1  @> \hat{\psi}^{7}_{0;1}>>   \pi_0 \D (D^{8})\\
 	@V P\otimes id VV       @V P VV\\
 	\mathrm{Coker} J_{8}\otimes \pi^s_1 @>precomposition >> 	\mathrm{Coker} J_{9}
 \end{CD}$$ induces the $2$-torsion element $\theta_2^6:=\psi^6_{1,1}(\theta_1^6\otimes \eta)\in \pi_2 \D (D^6) $.  
The following commutative diagram
\[
\begin{CD}
	\Z_{2} \tau_2^6 \oplus \Z_{2} \theta_2^6 @>\cong>> \Theta_{9}/bP_{10}\\
	@VVV       @A P AA\\
	\pi_{2} \D (D^6)  @>>> \pi_2 \bD(D^6)
\end{CD}
\] 
indicates that the homomorphism $ \pi_{2} \D (D^6)  \longrightarrow \pi_2 \bD(D^6)\cong \Theta_{9}\longrightarrow \Theta_9/bP_{10}$ is split surjective.
\end{proof}

The case of dimension 16 is similar to the case of dimension 8. Choose $\tilde{t}_7\in \pi_7SO(7)$ with $J_{7,7}(\tilde{t}_7)=\sigma'\in \pi_{14}S^7$. Let $\eta_7\in \pi_8 S^7$ be the generator and $\gamma_{7,8}\in \pi_7 SO(8)$ be the second generator
with $J_7 s_{\infty}\gamma_{7,8}=J_7\gamma_7=\sigma + \alpha_2+\alpha_{1,5}\in \pi^s_7$.
For the following Milnor--Novikov pairing
\[
\begin{CD}
	\pi_{8}SO(7) @.\otimes @.\pi_7 SO(8) @>\sigma_{8,7}^{(0,0)}>>  \pi_{1}\D (D^{14})\\
	@Vs_0VV         @.           @Vs_0VV                @VV \lambda_{1,1}^{15} V\\
	\pi_{8}SO(7)   @.\otimes @. \pi_7 SO(8) @>\sigma_{8,7}>>     \pi_0 \D (D^{15})\cong \Theta_{16},     
\end{CD}
\]
let $\theta_1^{14}=\sigma_{8,7}^{(0,0)}(\tilde{t}_7\eta_7,15\gamma_{7,8})\in \pi_{1}\D (D^{14})$, one has
\begin{prop}\label{dim14z2}
	\begin{itemize}
		\item[(a)] 	$\lambda_{1,1}^{15}(\theta_1^{14})=\sigma_{8,7}(\tilde{t}_7\eta_7,15\gamma_{7,8})=\Sigma^{16}_{[\eta^{*}]}\in \Theta_{16}$;
		\item[(b)] $\theta_1^{14}$ is of order 2;
		\item[(c)] the forgetful homomorphism $\pi_{1} \D (D^{14})  \longrightarrow \pi_1 \bD(D^{14})$ is split surjective.
			\end{itemize}

\end{prop}
\begin{proof}
Write $\alpha$ for $ \tilde{t}_7\eta_7$. Let $\beta=15\gamma_{7,8}$ to cancel the odd torsion. We see 
\begin{itemize}
	\item[(1)] $J_{7,7}(\alpha)=J_{7,7} \tilde{t}_7\eta_7=\sigma'\eta_{14}$ with $\Sigma^{\infty}\sigma'=2\sigma$ and $\Sigma^{\infty}\eta_{14}=\eta$;
	\item[(2)] $J_8(s_{\infty}\alpha)=\Sigma^{\infty}J_{7,7} \tilde{t}_7\eta_7=2\sigma\eta=0$ implies $s_{\infty}\alpha=0$;
		\item[(3)]  $\Sigma^{\infty}J_{7,8}\beta=15(\sigma + \alpha_2+\alpha_{1,5})=15\sigma$.
	%	\item[(4)] By Lemma \ref{relationhomotopy}, $\sigma\circ \pi_9^s=\Z_{2}\{\eta\sigma\overline{\nu} ,\eta\sigma \epsilon ,\sigma \mu \}=\Z_{2}\{\eta\rho\}=\mathrm{Im}J_{16}$.
\end{itemize}
By Theorem \ref{FrankSmith}, 
\[
P\sigma_{8,7}(\tilde{t}_7\eta_7,15\gamma_{7,8})\in - q(<15\sigma,2\sigma,\eta>)\subset [-15\eta^*]=[\eta^*].
\]
In dimension 16, the composition map $P\colon \Theta_{16}\longrightarrow \mathrm{Coker}J_{16}$ is an isomorphism and we have 
$$\lambda_{1,1}^{15}(\theta_1^{14})=\lambda_{1,1}^{15}(\sigma_{8,7}^{(0,0)}(\tilde{t}_7\eta_7,15\gamma_{7,8}))=\sigma_{8,7}(\tilde{t}_7\eta_7,15\gamma_{7,8})=\Sigma^{16}_{[\eta^{*}]}\in \Theta_{16}.$$ 
The equality $2\tilde{t}_7\eta_7=0$ implies $\theta_1^{14}$ is of order 2. For the subgroup generated by $ \theta_1^{14}\in \pi_{1} \D (D^{14})$, the following commutative diagram
\[
\begin{CD}
	\Z_{2} \theta_1^{14}@> \cong >> \Theta_{16}\\
	@VVV       @A \cong AA\\
	\pi_{1} \D (D^{14})  @>>> \pi_1 \bD(D^{14})
\end{CD}
\] 
implies that the forgetful homomorphism $\pi_{1} \D (D^{14})  \longrightarrow \pi_1 \bD(D^{14})$ is split surjective.

\end{proof}

\subsection{Dimension 10} This pairing comes from \cite[Example 2, Page 565]{frank1968invariant}. By Lemma \ref{unstableJ} (a), choose $\tilde{\alpha}_3$ with $J_{3,3}\tilde{\alpha_3}=\alpha_{1,S^3}$. For the following Milnor--Novikov pairing
\[
\begin{CD}
	\pi_6SO(3) @.\otimes @.\pi_3 SO(4) @>\sigma_{6,3}^{(0,2)}>>  \pi_{3}\D (D^{6})\\
	@Vs_0VV         @.           @Vs_2VV                @VV \lambda_{3,3}^{9} V\\
	\pi_6SO(3)   @.\otimes @. \pi_3 SO(6) @>\sigma_{6,3}>>     \pi_0 \D (D^{9})\cong \Theta_{10},     
\end{CD}
\]
let $\theta_3^6=\sigma_{6,3}^{(0,2)}(\tilde{\alpha}_3\alpha_{1,S^3}, 8\gamma_{3,4})\in \pi_3 \D (D^6)$, one has
\begin{prop}\label{dim10z3}
	\begin{itemize}
		\item[(a)] 	$\lambda_{3,3}^{9}(\theta_3^6)=\sigma_{6,3}(\tilde{\alpha}_3\alpha_{1,S^3}, 8\gamma_{3,6})=\Sigma^{10}_{[\beta_1]}\in \Theta_{10}$;
		\item[(b)] $\theta_3^6$ is of order 3;
		\item[(c)] the homomorphism $ \pi_{3} \D (D^6)  \longrightarrow \pi_3 \bD(D^6) $ is split surjective.
		\end{itemize}

\end{prop}
\begin{proof}
We write $\alpha$ for $\tilde{\alpha}_3\alpha_{1,S^3}$, and let $\beta=8\gamma_{3,4}$ to cancel the 2-torsion.
\begin{itemize}
	\item[(1)] $J_{6,3}\alpha=J_{6,3} \tilde{\alpha}_3\alpha_{1,S^3}=\alpha_{1,S^{3}}\alpha_{1,S^{6}}$, $\Sigma^{\infty}\alpha_{1,S^{3}}=\alpha_{1}$, $\Sigma^{\infty}\alpha_{1,S^{6}}=\alpha_{1}$;
	\item[(2)] $s_{\infty}\alpha=0\in \pi_6SO=0$;
	\item[(3)] $\Sigma^{\infty} J_{3,4}(\beta)=8(\nu+\alpha_{1})=2\alpha_{1}$.
%	\item[(4)] $\alpha_{1}\circ \pi^s_7=\{\alpha_{1}\alpha_2\}=0$ (e.g. \cite[Page 180]{toda1962composition}).
\end{itemize}
By Theorem \ref{FrankSmith}, one has
\[
P\sigma_{6,3}(\tilde{\alpha}_3\alpha_{1,S^3}, 8\gamma_{3,6})\in - q(<2\alpha_{1},\alpha_{1},\alpha_{1}>)\subset[-2\beta_1]=[\beta_1].
\]
In dimension 10, the composition map $P$ is an isomorphism and one has
$$\lambda_{3,3}^{9}(\theta_3^6)=\sigma_{6,3}^{(0,3)}(\tilde{\alpha}_3\alpha_{1,S^3}, 8\gamma_{3,4})=\sigma_{6,3}(\tilde{\alpha}_3\alpha_{1,S^3}, 8\gamma_{3,6})=\Sigma^{10}_{[\beta_1]}.$$
The equality $3\tilde{\alpha}_3\alpha_{1,S^3}=0$ implies that $\theta_3^6$ is of order 3. For the splitness of the homomorphism 
$$ \pi_{3} \D (D^6)  \longrightarrow \pi_3 \bD(D^6) ,$$
note that, by \cite[Theorem 1.1]{crowley2018harmonic}, there exists a
$2$-torsion element $\tau^6_3\in\pi_{3}\D (D^6) $ with $\alpha$-invariant 1. The subgroup $\Z_{2}\tau^6_3 \oplus \Z_{3} \theta_3^6\subset  \pi_{3}\D (D^6)$ fits into the following commutative diagram
\[
\begin{CD}
	\Z_{2}\tau^6_3 \oplus \Z_{3} \theta_3^6 @> \cong >>\Theta_{10}\\
	@VVV       @A \cong AA\\
	\pi_{3} \D (D^6)  @>>> \pi_3 \bD(D^6).
\end{CD}
\] 
Therefore the homomorphism $ \pi_{3} \D (D^6)  \longrightarrow \pi_3 \bD(D^6) $ is split surjective. 
\end{proof}
\begin{cor}\label{corodim13element}
 The homomorphism $ \pi_{6} \D (D^6)  \longrightarrow \pi_6 \bD(D^6) $ is split surjective.
\end{cor}
\begin{proof}
Let $\alpha_{1}\in \pi^s_3$ be the generator. By Proposition \ref{Compsitiondiagram}, the following composition diagram
$$	\begin{CD}
	\pi_3 \D (D^6) \otimes \pi^s_3 @> \psi^6_{3;3}>> 	\pi_{6} \D (D^6) \\
	@V \lambda^9_{3,3}\otimes id VV                                @V\lambda^{12}_{6,6} VV\\
	\pi_0 \D(D^{9})\otimes \pi^s_3  @> \hat{\psi}^{9}_{0;3}>>   \pi_0 \D (D^{12})\cong \Theta_{13}\\
	@V P\otimes id VV       @V P V \cong V\\
	\mathrm{Coker} J_{10}\otimes \pi^s_3 @>precomposition >> 	\mathrm{Coker} J_{13}
\end{CD}$$
induces a 3-torsion element $\theta_6^6:=\psi_{3,3}^6(\theta^6_3\otimes \alpha_{1})\in \pi_6 \D (D^6)$.  
The following commutative diagram
\[
\begin{CD}
	\Z_{3} \theta_6^6@> \cong >> \Theta_{13}\\
	@VVV       @A \cong AA\\
	\pi_{6} \D (D^6)  @>>> \pi_6 \bD(D^6)
\end{CD}
\] 
implies that the homomorphism $ \pi_{6} \D (D^6)  \longrightarrow \pi_6 \bD(D^6) $ is split surjective.
\end{proof}
\subsection{Dimension 15 and 18} The first pairing in dimension 15 comes from \cite[Example 3, Page 565]{frank1968invariant}. For the following Milnor--Novikov pairing
\[
\begin{CD}
	\pi_{11}SO(3) @.\otimes @.\pi_3 SO(4) @>\sigma_{11,3}^{(0,7)}>>  \pi_{8}\D (D^{6})\\
	@Vs_0VV         @.           @Vs_7VV                @VV \lambda_{8,8}^{14} V\\
	\pi_{11}SO(3)   @.\otimes @. \pi_3 SO(11) @>\sigma_{11,3}>>     \pi_0 \D (D^{14})\cong \Theta_{15},     
\end{CD}
\]
let $\theta_8^6=\sigma_{11,3}^{(0,7)}(\tilde{t}_3\epsilon_3,9\gamma_{3,4})\in \pi_8 \D (D^6)$, one has
\begin{prop}\label{MNpair15}

	\begin{itemize}
		\item[(a)] $P(\lambda_{8,8}^{14}(\theta_8^6))=P(\sigma_{11,3}(\tilde{t}_3\epsilon_3,9\gamma_{3,11}))=[\eta \kappa]\in \mathrm{Coker}J_{15}$.
		\item[(b)] $\theta_8^6$ is of order 2;
		\item[(c)] the composition $ \pi_{8} \D (D^6)  \longrightarrow \pi_8 \bD(D^6)\cong \Theta_{15}\longrightarrow \Theta_{15}/bP_{16}$ is split surjective. 
	\end{itemize}

\end{prop}
\begin{proof}
Note that  $J_{11,3} \tilde{t}_3\epsilon_3=\nu'\epsilon_6$ with $\Sigma^{\infty}\epsilon_3=\epsilon$. The equality $2\tilde{t}_3\epsilon_3=0$ yields $s_{\infty}(\tilde{t}_3\epsilon_3)=0$. The elements $\nu$, $2\nu$, and $\epsilon$ satisfy the Toda condition. By Theorem \ref{FrankSmith}, 
\[
P(\lambda_{8,8}^{14}(\theta_8^6))=	P(\lambda_{8,8}^{14}(\sigma_{11,3}^{(0,7)}(\tilde{t}_3\epsilon_3,9\gamma_{3,4}))=P(\sigma_{11,3}(\tilde{t}_3\epsilon_3,9\gamma_{3,11}))\in -q(<\nu, 2\nu, \epsilon>)=[-\eta\kappa]=[\eta\kappa].
\]
As $2\epsilon_3=0$, the order of $\theta_8^6$ is equal to 2.  
The following commutative diagram
\[
\begin{CD}
	\Z_{2} \theta_8^6 @>\cong>> \Theta_{15}/bP_{16}\\
	@VVV       @A P AA\\
	\pi_{8} \D (D^6)  @>>> \pi_8 \bD(D^6)
\end{CD}
\] 
implies that the composition $ \pi_{8} \D (D^6)  \longrightarrow \pi_8 \bD(D^6)\cong \Theta_{15}\longrightarrow \Theta_{15}/bP_{16}$ is split surjective.

\end{proof}

In dimension 18, for the following Milnor--Novikov pairing
\[
\begin{CD}
	\pi_{10}SO(7) @.\otimes @.\pi_7 SO(8) @>\sigma_{10,7}^{(0,2)}>>  \pi_{3}\D (D^{14})\\
	@Vs_0VV         @.           @Vs_2VV                @VV \lambda_{3,3}^{17} V\\
	\pi_{10}SO(7)   @.\otimes @. \pi_7 SO(10) @>\sigma_{10,7}>>     \pi_0 \D (D^{17})\cong \Theta_{18},     
\end{CD}
\]
let $\theta_3^{14}=\sigma_{10,7}^{(0,2)}(\tilde{t}_7\nu_7,15\gamma_{7,8})\in \pi_3 \D (D^{14})$, one has
\begin{prop}\label{dim18z2}

	\begin{itemize}
	\item[(a)] $\lambda_{3,3}^{17}(\theta_3^{14})=\sigma_{10,7}(\tilde{t}_7\nu_7,15\gamma_{7,10})=\Sigma^{18}_{[\nu^{*}]}\in \Theta_{18};$\\
	\item[(b)] $\theta_3^{14}$ is of order 8;
	\item[(c)] the homomorphism $ \pi_{3} \D (D^{14})  \longrightarrow \pi_3 \bD(D^{14})$ is split surjective.
	\end{itemize}
\end{prop}
\begin{proof}
For convenience, we choose $15\gamma_{7,8}$ to cancel the odd torsion. Note that  $J_{10,7}\tilde{t}_7\nu_7=\sigma'\nu_{14}$ and $s_{\infty}(\tilde{t}_7\nu_7)=0\in \pi_{10}SO=0$. The elements $\sigma,2\sigma$ and $\nu$ satisfy the Toda condition.
	By Theorem \ref{FrankSmith}, 
	\[
P(\lambda_{3,3}^{17}(\theta_3^{14}))=P(\lambda_{3,3}^{17}(\sigma_{10,7}^{(0,2)}(\tilde{t}_7\nu_7,15\gamma_{7,8}))=P(\sigma_{10,7}(\tilde{t}_7\nu_7,15\gamma_{7,10}))\in -q(<15\sigma, 2\sigma, \nu>)=[\nu^{*}].	\]
The composition map $P$ is an isomorphism in dimension 18, and (a) follows. The fact $8\nu_7=0$ implies the order of $\theta_3^{14}$ is at most 8, and the fact $8\nu^*=0$ implies the order of $\theta_3^{14}$ is at least 8. Hence, (b) follows.
Moreover, choose the 
$2$-torsion element $\tau_3^{14}\in \pi_{3} \D (D^{14})$ with $\alpha$-invariant 1 (\cite[Theorem 1.1]{crowley2018harmonic}), the following commutative diagram
\[
\begin{CD}
	\Z_{8} \theta_3^{14}\oplus \Z_{2} \tau_3^{14} @> \cong >> \Theta_{18}\\
	@VVV       @A \cong AA\\
	\pi_{3} \D (D^{14})  @>>> \pi_3 \bD(D^{14})
\end{CD}
\]
implies the homomorphism $ \pi_{3} \D (D^{14})  \longrightarrow \pi_3 \bD(D^{14})$ is split surjective.
\end{proof}

\section{Proof of Theorem \ref{Theorem A} and Theorem \ref{thmpi1}}\label{SectionProofMainthm}
In this section, we give the proof of our main results based on the methods in Section \ref{SectionMNTodaCompoistion} and the Milnor--Novikov pairings constructed in Section \ref{SectionMilnor}.
\subsection{Proof of Theorem \ref{Theorem A}}
\subsubsection{Proof of Theorem \ref{Theorem A} (a): $\Gamma_{4}^{10}$, $\Gamma_{7}^{13}$, and $\Gamma_{4}^{18}$}

\

\paragraph{$\bm{\Gamma_{4}^{10}}$} By Proposition \ref{dim10z3} (a), the 3-torsion element $\Sigma^{10}_{[\beta_1]}\in \Gamma_{4}^{10}$. According to \cite[Page 1097]{crowley2018harmonic}, the group $\Gamma_{4}^{10}$ also contains a 2-torsion element. It follows that $\Gamma_{4}^{10}=\Gamma_{2}^{10}\cong \Z_{2}\oplus \Z_{3}.$

\

\paragraph{$\bm{\Gamma_{7}^{13}}$} According to the proof of Corollary \ref{corodim13element}, 
for the 3-torsion element $\theta_6^6\in \pi_6 \D (D^6)$, we have
\[
P(\lambda^{12}_{6,6} (\theta_6^6))=P (\hat{\psi}^{9}_{0;3}(\Sigma^{10}_{[\beta_1]}\otimes\alpha_{1}))=[\alpha_{1}\beta_1]\in \mathrm{Coker} J_{13}.
\]
The composition map $P$ is an isomorphism in dimension 13 and we have $\Gamma_{7}^{13}=\Gamma_{2}^{13}\cong \Z_{3}$.

\

\paragraph{$\bm{\Gamma_{5}^{18}}$} By Proposition \ref{dim18z2}-(a), the element $\Sigma^{18}_{[\nu^{*}]}$ of order 8 belongs to $ \Gamma_{4}^{18}$. Furthermore, the fact that $\pi_{10}SO(6)\rightarrow \pi_{10}SO(7)$ is surjective implies $\Sigma^{18}_{[\nu^{*}]}\in \Gamma_{5}^{18}$. The 2-torsion element $\Sigma^{18}_{[\eta\overline{\mu}]}\in \Gamma_{11}^{18}$ (\cite[Page 1097]{crowley2018harmonic}). It follows that $\Gamma_{5}^{18}=\Gamma_{2}^{18}\cong \Z_{8}\oplus  \Z_{2}$.

\subsubsection{Proof of Theorem \ref{Theorem A} (b): $\Gamma_{3}^{9}$, $\Gamma_{9}^{15}$ and $\Gamma_{5}^{17}$} \label{proofofAb}
\paragraph{$\bm{\Gamma_{3}^{9}}$}
According to the proof of Corollary \ref{corodim9}-(a), the 2-torsion element $\theta_2^6$ admits 
\[
P(\lambda^8_{2,2}(\theta_2^6)) =P(\hat{\psi}^{8}_{0;1}(\Sigma^{8}_{[\epsilon]}\otimes\eta))=P(\epsilon)\circ \eta=[\epsilon\eta]\in \mathrm{Coker}J_{9}.
\]
As the group $\Gamma^9_3$ also contains a 2-torsion element with $\alpha$-invariant 1  (see \cite[Page 1097]{crowley2018harmonic} for the definition of the $\alpha$-invariant), we see $(\Z_{2})^2\subset \Gamma_{3}^9$.

\
\paragraph{$\bm{\Gamma_{9}^{15}}$}
By Proposition \ref{MNpair15}, the 2-torsion element $\theta_8^6 \in \pi_8 \D (D^6)$ satisfies 
$$P(\lambda_{8,8}^{14}(\theta_8^6))=[\eta\kappa]\in \coker J_{15} .$$ It follows that $\Z_{2}\subset \Gamma_{9}^{15}$.

\

\paragraph{$\bm{\Gamma_{5}^{17}}$}
%For $\Sigma^{16}_{[\eta^*]}\in \Gamma_{2}^{16}$, the 2-torsion element $\psi^{16}_{0;1}(\Sigma^{16}_{[\eta^*]}\otimes\eta)\in %\Gamma_{3}^{17}$ admits
%$$P (\psi^{16}_{0;1}(\Sigma^{16}_{[\eta^*]}\otimes\eta))=[\eta^*\eta].$$
First, for the unique 2-torsion element $\Sigma^{14}_{[\kappa]}\in \Gamma_{2}^{14}\cong \Theta_{14}$, choose the 2-torsion element $\theta_1^{12}\in \pi_{1}\D(D^{12})\cong \Theta_{14}$ with $\lambda_{1,1}^{12}(\theta_1^{12})=\Sigma^{14}_{[\kappa]}$ (\cite[Theorem 1.1]{wang2024short}). We can use the following commutative diagram
$$	\begin{CD}
	\pi_1 \D (D^{12}) \otimes \pi^s_3 @> \psi^{12}_{1;3}>> 	\pi_{4} \D (D^{12}) \\
	@V \lambda^{13}_{1,1}\otimes id VV                                @V\lambda^{16}_{4,4} VV\\
	\pi_0 \D(D^{13})\otimes \pi^s_3  @> \hat{\psi}^{13}_{0;3}>>   \pi_0 \D (D^{16})\cong \Theta_{17}\\
	@V P\otimes id VV       @V P V  V\\
	\mathrm{Coker} J_{14}\otimes \pi^s_3 @>precomposition >> 	\mathrm{Coker} J_{17}
\end{CD}$$
to construct the 2-torsion element $\theta_4^{12}:=\psi^{12}_{1;3}(\theta_1^{12}\otimes\nu)\in\pi_{4} \D (D^{12}) $ such that 
$$P(\lambda_{4,4}^{16}(\theta_4^{12}))=P (\hat{\psi}^{13}_{0;3}(\Sigma^{14}_{[\kappa]}\otimes\nu))=[\kappa\nu]\in \coker J_{17} .$$
Moreover, according to \cite[Lemma 2.2, Page 3702]{crowley2023derivative} and \cite[Page 1097]{crowley2018harmonic}, the homotopy sphere $\Sigma_{[\eta\eta^*]}\in \Gamma_{6}^{17}$ and the homotopy sphere $\Sigma_{[\overline{\mu}]}\in \Gamma_{6}^{17}$, we see $(\Z_{2})^3\subset \Gamma_{5}^{17}$.

\subsubsection{Proof of Theorem \ref{Theorem A} (c): $\Gamma_{4}^{4k+3}$ for $k\geq 3$}

%By \cite[Page 32]{ANTONELLI1972}, we see $\Z_{2032}\subset bP_{16}\cap \Gamma_{4}^{15}$. It follows that
%$\Z_{2}\lambda_{8,8}^{14}(\theta_8^6)\oplus \Z_{2032}\subset \Gamma_{4}^{15}\subset %\Z_{2}\lambda_{8,8}^{14}(\theta_8^6)\oplus \Z_{4064}\cong \Gamma_{3}^{15}$.

\

%\paragraph{$\Gamma_{4}^{4k+3}$ for $k\geq 3$}
The style of the proof of Theorem \ref{Theorem A}-(c) is different from that of (a) and (b). First, we need a lemma concerning the relative group $\frac{\widetilde{\mathrm{Diff}}_{\partial}(D^n)}{\D(D^n)}$.
\begin{lem}\label{pi3diffz2}
$\pi_3(\frac{\widetilde{\mathrm{Diff}}_{\partial}(D^n)}{\D(D^n)} ) \cong \Z_2$ when $n\geq 11$.
\end{lem}
\begin{proof}
The proof is similar to the proof of Proposition 2.2 in \cite{wang2024short}. 
Write $$\mathrm{C}(D^n):=\{f\in \mathrm{Diff}(D^n\times [0,1])  \ | \ f|_{D^n\times 0 \sqcup \partial D^n\times [0,1]}=id\}$$ 
for the concordance group of $D^n$. According to \cite[Page 915]{rognes2002two}, \cite[Page 174, 178-180]{rognes2003smooth} and \cite[Page 7]{igusa1988stability}, the second homotopy group $\pi_2 \mathrm{C}(D^n)$ of $\mathrm{C}(D^n)$
is zero when $n\geq 11$. For $p\geq 0$, the $E^1$-page of the Hatcher spectral sequence (see \cite[Page 4069]{wang2024short} for example) satisfies 
$$E^1_{p,0}=0, E^1_{p,1}\cong \Z_{2}, \text{and} \ E^1_{p,2}=0$$
with trivial differential $d^1\colon E^1_{1,1}\longrightarrow E^1_{0,1}$ (\cite[Lemma 2.1]{wang2024short}). It follows that 	$\pi_3(\frac{\widetilde{\mathrm{Diff}}_{\partial}(D^n)}{\D(D^n)} ) \cong \Z_2$ when $n\geq 11$.
\end{proof}

According to \cite[Theorem 1.4]{wang2024short}, the quotient groups  $\pi_2\widetilde{\mathrm{Diff}}_{\partial}(D^{4k})/\Gamma_3^{4k+3}\cong \Theta_{4k+3}/\Gamma_3^{4k+3}$ and 
$\pi_3\widetilde{\mathrm{Diff}}_{\partial}(D^{4k-1})/\Gamma_4^{4k+3}\cong \Theta_{4k+3}/\Gamma_4^{4k+3}$ in the short exact sequence (\ref{maineq}) are nontrivial subgroups of $\Z_2$. It follows that $\Gamma_{4}^{4k+3}=\Gamma_{3}^{4k+3}$ when $k\geq 3$.
\subsection{Proof of Theorem \ref{thmpi1}}

\subsubsection{Proof of Theorem \ref{thmpi1} (a)}\label{ProofthmB}
First, note that the case $n=10$ is trivial; the cases $n=6$ and $n=14$ have been proved in Proposition \ref{dim8z2}-(c) and Proposition \ref{dim14z2}-(c); and the case $n=12$ has been proved in \cite[Theorem 1.1]{wang2024short}.  We only need to consider the cases $n=7,8,11$, and $15$.

\

\paragraph{$\bm{\pi_1 \D (D^7)}$}By the definition of the Gromoll homomorphism (\cite[Section 1]{crowley2013gromoll}) and the proof of Corollary \ref{corodim9}, we have the following commutative diagram
\[
\xymatrix{
\Z_{2}\tau_2^6 \oplus \Z_{2} \theta_2^6 \ar@{^{(}->}[d] \ar@{^{(}->}[r] &   \pi_2 \D (D^6)  \ar[r]^{\lambda^{8}_{1,2}}\ar[d]  &   \pi_1 \D (D^7) \ar[d]\ar[dr] &\\
    \Theta_{9} \ar[r]^{\cong}                                  &   \pi_2 \bD (D^6)   \ar[r]^{\cong}  &  \pi_{1}\bD (D^7)   \ar[r]  & \Theta_9/bP_{10}.   
}
\]
It follows that the composition $\pi_1 \D (D^7)\longrightarrow \pi_1 \bD (D^7)\longrightarrow \Theta_9/bP_{10}$ is 
split surjective. 

On the other hand, by \cite[Theorem 1.4.3]{ANTONELLI1972} there exists a 2-torsion element $\kappa_1^7\in \pi_{1}\D (D^{7})$ which maps to the Kervaire sphere in $bP_{10}\cong \Z_{2}$. 

Denote $\lambda^8_{1,2}(\theta_2^6)$ and $\lambda^8_{1,2}(\tau_2^6)$ by $\theta_1^7$ and $\tau_1^7$, we have the following commutative diagram
\[
\begin{CD}
	\Z_{2}\kappa_{1}^7 \oplus \Z_{2} \theta_1^7\oplus\Z_{2}\tau_1^7  @> \cong >> \Theta_{9}\\
	@VVV       @A \cong AA\\
	\pi_{1} \D (D^{7})  @>>> \pi_1 \bD(D^{7}),
	\end{CD}
\] 
which shows that the homomorphism $\pi_{1} \D (D^{7})  \longrightarrow \pi_1 \bD(D^{7})$ is split surjective.

\

\paragraph{$\bm{{\pi_1 \D (D^{8})}}$} Similarly, according to the proof of Proposition \ref{dim10z3}-(c), we have the following commutative diagram
\[
\xymatrix{
	\Z_{2}\tau_3^6 \oplus \Z_{3} \theta_3^6 \ar@{^{(}->}[d] \ar@{^{(}->}[r] &   \pi_3 \D (D^6)  \ar[r]\ar[d]  &   \pi_1 \D (D^{8}) \ar[d]\\
	\Theta_{10} \ar[r]^{\cong}                                  &   \pi_3 \bD (D^6)   \ar[r]  & 	\pi_{1}\bD (D^{8}).    
}
\]
Write $\tau_1^8$ for $\lambda^9_{2,3}(\tau_3^6)$ and  $\theta_1^8$ for $\lambda^9_{2,3}(\theta_3^6)$,
it follows that the composition  $$\Z_{2}\tau_1^8 \oplus \Z_{3} \theta_1^8 \subset \pi_1 \D (D^{8})\longrightarrow \pi_1 \bD (D^{8})\cong\Theta_{10}$$ is 
split surjective. 

\

\paragraph{$\bm{\pi_1 \D (D^{11})}$} By Corollary \ref{corodim13element}, we have the following commutative diagram
\[
\xymatrix{
	\Z_{3} \theta_6^6 \ar[d]^{\cong} \ar@{^{(}->}[r] &   \pi_6 \D (D^6)  \ar[r]\ar[d]  &   \pi_1 \D (D^{11}) \ar[d]\\
	\Theta_{13} \ar[r]^{\cong}                                  &   \pi_6 \bD (D^6)   \ar[r]  &  \pi_{1}\bD (D^{11}).    
}
\]
Write $\theta_1^{11}$ for $\lambda^{12}_{5,6}(\theta_6^6)$, it follows that the composition $\Z_{3}\theta_{1}^{11}\subset\pi_1 \D (D^{11})\longrightarrow \pi_1 \bD (D^{11})\cong\Theta_{13}$ is 
split surjective. 

\

\paragraph{$\bm{\pi_1 \D (D^{15})}$} 
First, by \cite[Section 1.4.1]{ANTONELLI1972}, the plumbing construction of two tangent bundles of $S^9$
produces a 2-torsion element $\kappa_{1}^{15}\in \pi_1 \D (D^{15})$ which maps to the Kervaire sphere in $\Gamma_{2}^{17}$.

Let $\tau_{10}^6\in \pi_{10} \D (D^6)$ be the 2-torsion element with $\alpha$-invariant 1 (\cite[Theorem 1.1]{crowley2018harmonic} and write the 2-torsion element $\tau_1^{15}\in \pi_1 \D(D^{15})$ for the image of $\lambda^{16}_{9,10}$. 

In the proof of Theorem \ref{Theorem A}-(b) (See Section \ref{proofofAb}), we construct the 2-torsion element $\theta_4^{12}\in \pi_4\D(D^{12})$ and write $\theta_1^{15}=\lambda^{16}_{3,4}(\theta_4^{12})\in \pi_1 \D(D^{15})$.

Similarly, we can construct the 2-torsion element  $\eta_2^{14}=\psi^{14}_{1,1}(\theta_1^{14}\otimes \eta)\in \pi_2 \D (D^{14})$ where $\theta_1^{14}\in \pi_{1}\D(D^{14})$ satisfies $\lambda^{15}_{1,1}(\theta_{1}^{14})=\Sigma_{[\eta^*]}$. We write $\eta_1^{15}$ for $\lambda^{16}_{1,2}(\eta_2^{14})\in \pi_1 \D(D^{15}) $.
Then we have the following commutative diagram
\[
\begin{CD}
	\Z_{2}\kappa_1^{15} \oplus \Z_{2} \tau_{1}^{15} \oplus \Z_{2} \eta_1^{15}\oplus  \Z_{2}\theta_{1}^{15}   @> \lambda^{16}_{1,1} >>\Theta_{17}\\
	@VVV       @A \cong AA\\
	\pi_{1} \D (D^{15})  @>>> \pi_1 \bD(D^{15}).
\end{CD}
\] 
It follows that $\pi_{1} \D (D^{15}) \longrightarrow \pi_1 \bD(D^{15})$ is split surjective.
\subsubsection{Proof of Theorem \ref{thmpi1} (b)}\label{proofThmBb}
By Proposition \ref{MNpair15}, the following commutative diagram holds
\[
\begin{CD}
	\Z_{2} \theta_8^6  @>\lambda^{14}_{7,8}>> \Theta_{15}   @> P >> \Theta_{15}/bP_{16}\\
	@VVV                             @A AA\\
	\pi_{8} \D (D^6)    @>\lambda^{14}_{7,8}>> \pi_{1} \D(D^{13}).                  @. 
\end{CD}
\] 
The composition $P\lambda^{14}_{7,8}\colon 	\Z_{2} \theta_8^6 \longrightarrow \Theta_{15}/bP_{16} $ is also an isomorphism. Write $\theta_{1}^{13}$ for $\lambda^{14}_{7,8}(\theta_8^6)$, it follows that the composition $\Z_{2} \theta_{1}^{13}\subset \pi_{1} \D(D^{13}) \longrightarrow \Theta_{15}/bP_{16}$ is split surjective.
\subsection{Proof of Corollaries}
\subsubsection{Proof of Corollary \ref{CorollaryC}} \label{ProofofCoroC}
For the group $\pi_1 \D (D^{15})$, consider the following short exact sequence,
\[
0\rightarrow \Theta_{18}/\Gamma^{18}_{3}\longrightarrow \pi_{2}\frac{\widetilde{\mathrm{Diff}}_{\partial}(D^{15})}{\D(D^{15})}\longrightarrow \pi_{1}\D(D^{15})\longrightarrow
\Theta_{17}\longrightarrow 0.
\]
By Theorem \ref{Theorem A}-(a) and Theorem \ref{thmpi1}-(a), we see $\Theta_{18}/\Gamma^{18}_{3}=0$ and the homomorphism $$\pi_{1}\D(D^{15})\longrightarrow
\Theta_{17}$$ is split surjective. It follows that 
\[
\Z_{2}  \oplus \Z_{2}\kappa_1^{15} \oplus \Z_{2} \tau_1^{15}\oplus \Z_{2}\theta_{1}^{15} \oplus \Z_{2} \eta_{1}^{15}\cong \pi_{1} \D (D^{15}).
\]

%The cases when $n=6,8,11,13$ and $14$ have been proved in Corollary \ref{Coropi1}.
%The cases when $n=10$ and $4k, k\geq 3$ have been proved in \cite[Theorem 1.1 \& Example 2.3]{wang2024short}.

%\

For the group $\pi_1 \D (D^{7})$,  consider the following short exact sequence
\[
0\longrightarrow \Theta_{10}/\Gamma^{10}_{3}\longrightarrow \pi_{2}\frac{\widetilde{\mathrm{Diff}}_{\partial}(D^7)}{\D(D^7)}\longrightarrow \pi_{1}\D(D^7)\longrightarrow
\Theta_{9}\longrightarrow 0.
\]
By Theorem \ref{Theorem A}-(a) and Theorem \ref{thmpi1}-(a), we see $\Theta_{10}/\Gamma^{10}_{3}=0$ and this sequence splits. It follows that $$\pi_{1} \D (D^7) \cong  \pi_2(\frac{\widetilde{\mathrm{Diff}}_{\partial}(D^7)}{\D(D^7)} )\oplus \Z_{2}\kappa_{1}^7 \oplus \Z_{2} \theta_1^7\oplus\Z_{2}\tau_1^7.$$

\

For the group {$\pi_1 \D (D^{9})$}, the short exact sequence becomes
\[
0\cong \Theta_{12}/\Gamma^{12}_{3}\longrightarrow \pi_{2}\frac{\widetilde{\mathrm{Diff}}_{\partial}(D^9)}{\D(D^9)}\longrightarrow \pi_{1}\D(D^9)\longrightarrow
\Theta_{11}\longrightarrow 0.
\]
Note that $\pi_{2}\frac{\widetilde{\mathrm{Diff}}_{\partial}(D^9)}{\D(D^9)}\cong \Z_{2}$ (\cite[Proposition 2.2]{wang2024short}) and the result follows.

\

\subsubsection{Proof of Corollary \ref{CorollaryD}}
When $n=4k-1\geq 15$, by Lemma \ref{pi3diffz2} and Theorem \ref{Theorem A}-(c), we have exact sequence
\[
0\longrightarrow \Z_{2}\longrightarrow \Z_2 \longrightarrow \pi_{2}\D(D^{4k-1})\longrightarrow \Gamma_3^{4k+2}\longrightarrow 0. 
\]
It follows that $\pi_{2}\D D^{4k-1}\cong \Gamma_3^{4k+2}$. Moreover, $\pi_{2}\D( D^{15})\cong \Gamma_3^{18}\cong \Theta_{18}$.

\subsubsection{Proof of Corollary \ref{CorollaryE}}
Let $M^n$ be a compact smooth manifold. Consider the homomorphism $\D (D^n)\hookrightarrow \D (M^n)$ defined by extension by the identity.
We have the following commutative diagram
\[
\begin{CD}
	\pi_3(\frac{\widetilde{\mathrm{Diff}}_{\partial}(D^{4k-1})}{\D(D^{4k-1})} ) @> trivial >> \pi_2\D(D^{4k-1}) @>>>   \pi_2\widetilde{\mathrm{Diff}}_{\partial}(D^{4k-1}) @>>>	\pi_2(\frac{\widetilde{\mathrm{Diff}}_{\partial}(D^{4k-1})}{\D(D^{4k-1})} )\\
	@V\cong VV   @VVV    @VVV  @VVV\\
	\pi_3(\frac{\widetilde{\mathrm{Diff}}_{\partial}(M^{4k-1})}{\D(M^{4k-1})} ) @>>> \pi_2\D(M^{4k-1}) @>>>   \pi_2\widetilde{\mathrm{Diff}}_{\partial}(M^{4k-1}) @>>>	\pi_2(\frac{\widetilde{\mathrm{Diff}}_{\partial}(M^{4k-1})}{\D(M^{4k-1})} ).\\
\end{CD}
\]
When $M^{4k-1}$ is 3-connected, by Morlet disjunction \cite[Corollary 3.2]{BurgheleaLashofRothenberg}, $$\pi_i(\frac{\widetilde{\mathrm{Diff}}_{\partial}(D^{4k})}{\D(D^{4k})} )\cong \pi_i(\frac{\widetilde{\mathrm{Diff}}_{\partial}(M^{4k})}{\D(M^{4k})} )$$ 
for $1\leq i \leq 3$, which implies the results by diagram chasing.

 \section{Proof of the Frank--Kosinski--Smith theorem}\label{SectionProofofFrankSmith}
 In this section, we give the detailed proof of the Frank--Kosinski--Smith theorem\footnote{We refer to Remark \ref{remarkofFKSthm} for the explanation of the name of ``the Frank--Kosinski--Smith theorem''.} (Theorem \ref{FrankSmith}). We follow Frank's ideas in \cite{frank1968invariant}, with techniques mainly from \cite{BrumfielI} and \cite{frank1967reducible}.
\subsection{Topological $(D^k,0)$-bundles and topological $J$-homomorphisms}
Let $M$ be a manifold. Denote by $\mathrm{Top}^{+}(M)$ the topological group of orientation-preserving homeomorphisms of $M$ with compact-open topology and denote by $\mathrm{Top}^{+}_A(M)$ the subgroup of $\mathrm{Top}^{+}(M)$ which sends the subspace $A\subset M$ onto itself. 

When $M=D^k$ and $A$ is the origin point of $D^k$, the topological group $\mathrm{Top}^{+}_{0}(D^k)$ is the structure group of $(D^k,0)$-bundles, i.e. topological $D^k$-bundles with zero-sections. Its classifying space $\mathrm{B}\mathrm{Top}^{+}_{0}(D^k)$ classifies oriented $(D^k,0)$-bundles.

Restricting to the interior points of $D^k$ and the boundary of $D^k$, we have two restriction homomorphisms $int\colon\mathrm{Top}^{+}_{0}(D^k)\longrightarrow \mathrm{Top}^{+}_{0}(\mathrm{int} (D^k))$ and $bd\colon\mathrm{Top}^{+}_{0}(D^k)\longrightarrow \mathrm{Top}^{+}(S^{k-1})$. 
Since $\mathrm{int}(D^k)$ is homeomorphic to $\R^k$ with the origin fixed, we can identify $\mathrm{Top}^{+}_{0}(\mathrm{int} (D^k))$ with $\mathrm{Top}^{+}_{0}(\R^k)$.

\begin{lem} The homomorphism $\pi_i\T_0 (D^k)\longrightarrow \pi_i\T_0 (\R^k)$ induced by the restriction homomorphism $int\colon\mathrm{Top}^{+}_{0}(D^k)\longrightarrow \mathrm{Top}^{+}_{0}(\mathrm{int} (D^k))\cong \mathrm{Top}^{+}_{0}(\R^k)$ is an isomorphism when $k>>i$.
%	\begin{itemize}
%		\item[(1)] $\mathrm{Top}^+(\R^k)\cong  \mathrm{Top}^{+}_{pt}(S^k)$.
%		\item[(2)]
%	\end{itemize}
\end{lem}
\begin{proof}
%	First we have commutative diagram
%	\[
%	\begin{CD}
%		\mathrm{Top}_0 D^k @> int >> \mathrm{Top}_0 \R^k\\
%		@V\simeq VV                                @V\simeq VV\\
%			\mathrm{Top} D^k @> int >> \mathrm{Top} \R^k
%	\end{CD}
%	\]
	First recall that the stabilization map $\pi_i\T_0 \R^k\longrightarrow \pi_i\T_0 \R^{k+1}$ is an isomorphism when $k>>i$ (e.g. \cite[Page 37]{BurgheleaLashof}). The one point compactification of $\R^k$ induces the group isomorphism 
	$$\mathrm{Top}^{+}(\R^k)\cong \mathrm{Top}^{+}_{pt}(S^k).$$
	The homotopy fibration $\T_{pt}(S^k) \longrightarrow \T(S^k) \longrightarrow S^k$ induces isomorphism 
	$\pi_i \T_{pt} S^k \cong \pi_i \T S^k $ when $k>>i$.  
Then the following commutative diagram
	\[
	\begin{CD}
	\T_0(\R^k) @>\simeq>>		\T(\R^k) @>\cong >> \T_{pt}(S^k)  @>>> \T(S^k)\\
	@V\times \R VV	@V\times \R VV    @V \Sigma VV            @V\Sigma VV\\
	\T_0(\R^{k+1}) @>\simeq >>		\T(\R^{k+1}) @>\cong >> \T_{pt}(S^{k+1})  @= \T_{pt}(S^{k+1}) \\
	\end{CD}
	\]
	indicates that the homomorphism $ \pi_i \T S^k \longrightarrow \pi_i \T_{pt} S^{k+1}$ induced by the suspension map $\Sigma$ is isomorphic when $k>>i$.
	
	Second,  note that the inclusion $\T_0 \R^k\hookrightarrow \T\R^k$ is a homotopy equivalence.  The Alexander trick \cite[Page 227]{browder1966open} implies the radical extension $\T S^{k-1}\rightarrow \T_0 D^k $ and the inclusion 	$ \T_0 D^k \hookrightarrow \T D^k$ are also homotopy equivalences. Therefore the following commutative diagram
	\[
	\xymatrix{
		\T S^{k-1} \ar[d]^{\simeq} \ar[rrd]^{\Sigma} &       \\
			\T_0 D^k  \ar[d]^{\simeq} \ar[r]^{int}  & \T_0 \R^k \ar[d]^{\simeq}  \ar[r]^{\simeq} & \T_{pt} S^k \ar@{=}[d]\\
			\T D^k \ar[r]^{int}    & \T \R^k \ar[r]^{\cong}  & \T_{pt} S^k \\ }
			\]
implies $\pi_i \mathrm{Top}^{+}_{0}(D^k)\cong \pi_i \mathrm{Top}^{+}_{0} (\R^k)$ when $k>>i$.
	\end{proof}

Since
$\mathrm{Top}^{+}(\R^k)\cong \mathrm{Top}^{+}_{pt}(S^k)$, one has
natural inclusions
\[
\mathrm{Top}^{+}_{0}(D^k)\longrightarrow \mathrm{Top}^{+}_{0}(\R^k) \longrightarrow \mathrm{Top}^{+}(\R^k)\cong
\mathrm{Top}_{pt}(S^k)\longrightarrow \Omega^k S^k,
\]
which induce the topological (unstable) $J$-homomorphism
\[
J^{\mathrm{T}}_{n}\colon \pi_n \mathrm{Top}^{+}_{0}(D^k) \longrightarrow \pi_{n+k} S^k.
\]

\subsection{Attaching maps of Thom spaces of certain topological $(D^k,0)$-bundles}
Let $\xi_{D^k}$ be an oriented $(D^k,0)$-bundle over a finite CW-complex $X$ with structure group $\mathrm{Top}^{+}_{0}(D^k)$. The internal points $\mathrm{int}(\xi_{D^k})$ of $\xi_{D^k}$ is a $(\mathrm{int}(D^k),0)$-bundle over $X$ with structure group $\mathrm{Top}^{+}_{0}(\mathrm{int} (D^k))$, which can also be identified with a $(\R^k,0)$-bundle with structure group $\mathrm{Top}^{+}_{0}(\R^k)$.  
The boundary points $\xi_{D^k}-\mathrm{int}(\xi_{D^k})$ becomes an oriented $S^{k-1}$-bundle with structure group $\mathrm{Top}^{+}(S^{k-1})$. 
\begin{defn}
	The Thom space $T(\xi_{D^k})$ of the topological $(D^k,0)$-bundle $\xi_{D^k}$ is defined to be the quotient space $\xi_{D^k}/(\xi_{D^k}-\mathrm{int}(\xi_{D^k})).$ 
\end{defn}

When $Y$ is a connected CW-complex whose dimension is at most $n-1$, let $X=Y\cup_{d}D^n$ with attaching map
$d\colon S^{n-1}\longrightarrow Y$. Let $\xi_{\alpha}$ and $\xi_{\beta}$ be two oriented topological $(D^k,0)$-bundles over $X$ with  $\xi_{\alpha}|_Y\cong \xi_{\beta}|_Y$. Denote the Thom spaces of $\xi_{\alpha}$ and $\xi_{\beta}$ by $T(\xi_{\alpha})$ and $T(\xi_{\beta})$, then $$T(\xi_{\alpha})=T(\xi_{\alpha}|_Y)\cup_{\phi_{\alpha}} D^{n+k}, \ T(\xi_{\beta})=T(\xi_{\beta}|_Y)\cup_{\phi_{\beta}} D^{n+k} $$ with attaching maps $[\phi_{\alpha}],[\phi_{\beta}]\in \pi_{n+k-1}T(\xi_{\alpha}|_Y)$. 

\begin{lem}[{\cite[Proposition 1]{frank1968invariant}}]\label{FranklemThomSpace}
 Suppose the attaching map $\phi_{\beta}$ of $T(\xi_{\beta})$ is null-homotopic, then $[\phi_{\alpha}]\in \mathrm{Image}(i_{\ast})$, where the homomorphism $$i_{\ast}\colon\pi_{n+k-1}S^k\longrightarrow \pi_{n+k-1}T(\xi_{\alpha}|_Y)$$ is induced by the inclusion $i\colon S^k\cong T(\xi_{\alpha}|_{pt})\hookrightarrow T(\xi_{\alpha}|_Y)$.
\end{lem}
\begin{proof}
We follow the proof in \cite[Page 14-18]{frank1967reducible}. Since $X=Y\cup_{d} D^n$, the bundle $\xi_{\alpha}=\xi_{\alpha}|_Y \cup \xi_{\alpha}|_{D^n}$ is equivalent to $\xi_{\alpha}|_Y \cup_{\tilde{\phi}_{\alpha}}
	D^n\times D^k $ with a fiber-preserving attaching map
	\[
	\begin{CD}
		\tilde{\phi}_{\alpha}\colon	S^{n-1}\times D^k@>>> \xi_{\alpha}|_Y\\
		@VVV            @VVV\\
		S^{n-1} @>d>>  Y. 
	\end{CD}
	\]
	As a quotient space of $\xi_{\alpha}$, the Thom space $T(\xi_{\alpha})$ is homeomorphic to $T(\xi_{\alpha}|_Y)\cup_{\phi_{\alpha}} D^{n+k}$ with the attaching map $\phi_{\alpha}\colon S^{n+k-1}=\partial (D^k\times D^n)=S^{k-1}\times D^n\cup D^{k}\times S^{n-1}\longrightarrow \xi_{\alpha}|_Y$ described as follows
	\[
	\phi_{\alpha}(x)=\left\{
	\begin{array}{cc}
		\text{base point of  $\xi_{\alpha}|_Y$}	&  x \in S^{k-1}\times D^n  \\
		c\circ\tilde{\phi}_{\alpha}	&  x\in  D^{k}\times S^{n-1}\\
	\end{array}
	\right.
	\]
	where $c$ denotes the quotient map $\xi_{\alpha}|_Y\longrightarrow T(\xi_{\alpha})|_Y$.
	
	For the bundles $\xi_{\alpha}=\xi_{\alpha}|_Y \cup_{\tilde{\phi}_{\alpha}}
	D^n\times D^k$ and $\xi_{\beta}=\xi_{\beta}|_Y \cup_{\tilde{\phi}_{\beta}}
	D^n\times D^k$, the bundle isomorphism  $\xi_{\alpha}|_{Y}\cong\xi_{\beta}|_{Y}$ induces
	the following fiber-preserving equivalence
	\[
	\begin{CD}
		\tilde{\Theta}_{(\xi_{\alpha},\xi_{\beta})}\colon	S^{n-1}\times D^k  @>>> 	S^{n-1}\times D^k\\
		@VVV            @VVV\\
		S^{n-1} @>>>  S^{n-1}, 
	\end{CD}
	\]
	with $\tilde{\Theta}_{(\xi_{\alpha},\xi_{\beta})}(x,y)=(x, \tilde{\phi}^{-1}_{\beta,x}\tilde{\phi}_{\alpha,x}(y))$, where $\tilde{\phi}_{\alpha,x}\colon D^k\longrightarrow \xi_{\alpha}|_{d(x)}$ is the origin-preserving homeomorphism 
	of $D^k$ onto the fiber of $d(x)$. Define 
	$\Theta_{\xi_{\alpha},\xi_{\beta}}\colon S^{n-1}\rightarrow \mathrm{Top}^{+}_0(D^k)$ to be
	\[
	x\mapsto \mathrm{pr}_2 \tilde{\Theta}_{(\xi_{\alpha},\xi_{\beta})}(x,-).
	\]

	%We also denote $\Theta_{(\xi_{\alpha},\xi_{\beta})}
	
	For the Thom spaces $T(\xi_{\alpha})=T(\xi_{\alpha}|_Y)\cup_{\phi_{\alpha}} D^{n+k}$ and $T(\xi_{\beta})=T(\xi_{\beta}|_Y)\cup_{\phi_{\beta}} D^{n+k}$, let $$i_{\ast}\colon\pi_{n+k-1}S^k\longrightarrow \pi_{n+k-1}T(\xi_{\alpha}|_Y)$$
	 be the homomorphism induced by the inclusion $i\colon S^k\cong T(\xi_{\alpha}|_{pt})\hookrightarrow T(\xi_{\alpha}|_Y)$. According to \cite[Page 17]{frank1967reducible} , one has
	\[
	[\phi_{\alpha}]-[\phi_{\beta}]=i_*J_{n-1}^{\mathrm{T}}(\Theta_{(\xi_{\alpha},\xi_{\beta})})\in \pi_{n+k-1}T(\xi_{\alpha}|_Y),
	\]
	where $J_{n-1}^{\mathrm{T}}(\Theta_{(\xi_{\alpha},\xi_{\beta})})$ is the (topological) $J$-homomorphism of $\Theta_{(\xi_{\alpha},\xi_{\beta})}$.

	When the attaching map $\phi_{\beta}$ is null-homotopic, which is equivalent to that $T(\xi_{\beta})$ is reducible \cite[Page 12-13]{frank1967reducible}, it follows that $[\phi_{\alpha}]=i_*J_{n-1}^{\mathrm{T}}(\Theta_{(\xi_{\alpha},\xi_{\beta})})\in \mathrm{Image}(i_{\ast})\subset\pi_{n+k-1}T(\xi_{\alpha}|_Y)$.
\end{proof}
	\begin{rem}\label{remmarkfibrepreserving}
	In Lemma \ref{FranklemThomSpace}, the homotopy class of the attaching map $\phi_{\alpha}$ depends on the choice of a fiber-preserving equivalence $\tilde{\phi}_{\alpha}$. If we choose another fiber-preserving equivalence $\tilde{\phi}_{\alpha'}$, let $\xi_{\alpha'}=\xi_{\alpha}|_Y \cup_{\tilde{\phi}_{\alpha'}} 	D^n\times D^k$ and we also have
	$$[\phi_{\alpha}]-[\phi_{\alpha'}]=i_*J_{n-1}^{\mathrm{T}}(\Theta_{(\xi_{\alpha},\xi_{\alpha'})}).$$
\end{rem}
\begin{ex}\label{ExampleTSn}
	When $Y=pt$, $X=Y\cup D^n\cong S^n$, and $\xi$ is a topological $(D^k,0)$-bundle over $X$ with the clutching map $f_{\xi}\colon S^{n-1}\longrightarrow \T_0 (D^k)$. 
	We have homotopy equivalence 
	$T(\xi)\simeq S^k\cup_{\phi_{\xi}} D^{n+k}$
	with $\phi_{\xi}=J_{n-1}^{\mathrm{T}}(f_{\xi})$. When $\xi$ is the disc bundle of an oriented $k$-dimensional vector bundle over $X$ with the clutching map $f_{\xi}\colon S^{n-1}\longrightarrow SO(k)$, one has $T(\xi)=S^k\cup_{\phi_{\xi}} D^{n+k}$ with $\phi_{\xi}\simeq J_{n-1,k}(f_{\xi})$, which coincides with Lemma 10.1 in \cite{ADAMS196621}.
\end{ex}

\begin{lem}\label{attachThomspace}
Let $Y=S^n$, $X=Y\cup_d D^{n+l}$ with $d=\Sigma g$ for $g\colon S^{n+l-2}\longrightarrow S^{n-1}$. Let $\xi$ be a vector bundle over $X$ and let $f_{\xi_Y}\colon S^{n-1}\longrightarrow SO(k)$ be the clutching map of the restriction bundle $\xi_Y$ of $\xi$ to $Y$. For any attaching map of the Thom space $T(\xi)=T(\xi_Y)\cup_{\phi_{\xi}} D^{n+l+k}$, one has
\[
c\circ \phi_{\xi}\simeq \Sigma^{k+1}g,
\]
where $c:T(\xi_Y)\longrightarrow T(\xi_Y)/S^k\cong S^{n+k}$ denotes the quotient map.
\end{lem}
\begin{proof}
The cofiber sequence
\[
\begin{CD}
	S^{n+l-1} @>d>>  Y@>>>  X
\end{CD}
\]
induces bundle homomorphism
\[
\begin{CD}
	S^{n+l-1}\times \R^k\cong d^* \xi_Y  @>>>  \xi_Y   @>>>  \xi   \\
	  @VVV    @VVV   @VVV\\
	  	S^{n+l-1} @>d=\Sigma g>>  Y@>>>  X.
\end{CD}
\]
Note that the clutching maps of $\xi_Y$ and $d^*\xi_Y$ are $f_{\xi_Y}$ and $f_{\xi_Y}g$.
According to \ref{ExampleTSn}, $T(\xi_Y)\simeq S^k\cup_{J_{n-1,k}(f_{\xi_Y})} D^{n+k}$ and $T(d^* \xi_Y)\simeq  S^k\cup_{J_{n+l-2,k}(f_{\xi_Y}g)} D^{n+l+k-1}\simeq S^k\cup_{J_{n-1,k}(f_{\xi_Y})(\Sigma^{k}g)} D^{n+l+k-1}$. The map 
$T(d)\colon T(d^*\xi_Y)\longrightarrow T(\xi_Y)$ admits the following homotopy commutative diagram
\[
\begin{CD}
S^{n+l+k-2} @>J_{n+l-2,k}( f_{\xi_Y}g)>>  S^k @>i>> T(d^* \xi_Y)\simeq  S^k\cup_{(\Sigma^{k}g)J_{n+l-2,k}(f_{\xi_Y})} D^{n+l+k-1} @>c>>  S^{n+l+k-1}\\
@VV \Sigma^k gV      @|           @V T(d) VV   @V\Sigma^{k+1}g VV\\     
S^{n+k-1}  @>J_{n-1,k}(f_{\xi_Y})>>  S^k  @>i>>  T(\xi_Y)\simeq S^k\cup_{J_{n-1,k}(f_{\xi_Y})} D^{n+k}  @>c>>  S^{n+k}. 
\end{CD}
\]
The attaching map $\phi_{\xi}$ of $T(\xi)=T(\xi_Y)\cup_{\phi_{\xi}} D^{n+l+k}$ admits the following homotopy commutative diagram
\[
\xymatrix{
\phi_{\xi}\colon S^{n+l+k-1} \ar[rrd]^{\rho}  \ar[r] &      S^{n+l+k-1}/(S^{k-1}\times D^{k-1})\ar[r]^-{\cong} &T(d^*\xi_Y)  \ar[r]^{T(d)}\ar[d]^{c} &   T(\xi_Y)\ar[d]^c \\
              &     &     S^{n+l+k-1}    \ar[r]^{\Sigma^{k+1} g}          &  S^{n+k}.\\}\]
Observe that deg$\rho$=1 and $c\phi_{\xi}=(\Sigma^{k+1}g)\rho\simeq \Sigma^{k+1}g$. 
\end{proof}

\begin{lem}[Frank, {\cite[Example 2, Page 24]{frank1967reducible}}]\label{TodaLemma}
Let $Y=S^n$, $X=Y\cup_{(\Sigma d_1)  d_2}D^{n+l_1+l_2}$ and  $X_1=Y\cup_{\Sigma d_1}D^{n+l_1}$ with maps
$$d_1\colon S^{n+l_1-2}\longrightarrow S^{n-1}, d_2\colon S^{n+l_1+l_2-1}\longrightarrow S^{n+l_1-1}.$$
Let $F\colon X\longrightarrow X_1$ be the map with $F|_Y=id_Y$ and $F|_{D^{n+l_1+l_2}}$ is the radical extension of $d_2$
\[
\begin{CD}
      S^{n+l_1+l_2-1}         @>>>  Y  @>>> X= Y\cup_{(\Sigma d_1)  d_2}D^{n+l_1+l_2}@>>> S^{n+l_1+l_2}\\
      @V d_2 VV        @|    @V F VV     @V \Sigma d_2 VV\\
       S^{n+l_1-1}         @>>>  Y  @>>> X_1= Y\cup_{(\Sigma d_1) }D^{n+l_1} @>>>  S^{n+l_1} .\\
\end{CD}
\]
Let $\xi$ be a $k$-dimensional vector bundle over $X_1$ whose clutching map of $\xi_Y$ over $Y$ is $f_{\xi_Y}$ and let $F^*\xi$ be the pull-back bundle over $X$
\[
\begin{CD}
\xi_Y @>>>	F^*\xi @>>> \xi\\
@VVV	@VVV   @VVV\\
Y	@>>>X  @>F>> X_1.
	\end{CD}
\]
For the attaching map $\phi_{
	F^*\xi}$ of $T(F^*\xi)=T(\xi_Y)\cup_{\phi_{F^*\xi}} D^{n+l_1+l_2+k}$, 
assume $k>>n+l_1+l_2$ and $[\phi_{F^*\xi}]\in \mathrm{Image} \ i_*\colon\pi_{n+l_1+l_2+k-1}S^k\longrightarrow  \pi_{n+l_1+l_2+k-1}T(\xi_Y)$ where $i_*$ is induced by the inclusion $i\colon S^k\longrightarrow T(\xi_Y)\simeq S^k\cup_{J_{n-1.k}(f_{\xi_Y})} D^{n+k}$, one has
\[
[\phi_{F^*\xi}]\in -i_* \{J_{n-1,k}(f_{\xi_Y}),\Sigma ^k d_1,\Sigma^{k-1}d_2 \},
\] 
where $\{J_{n-1,k}(f_{\xi_Y}),\Sigma ^k d_1,\Sigma^{k-1}d_2 \}$ is the Toda bracket of $J_{n-1,k}(f_{\xi_Y})$, $\Sigma ^k d_1$, and $\Sigma^{k-1}d_2$.
\end{lem}
\begin{proof}
%Let $X_1=Y\cup_{\Sigma d_1}D^{n+l_1}$ and let $F\colon X\longrightarrow X_1$ be the map with $F|_Y=id_Y$ and %$F|_{D^{n+l_1}}$ is the radical extension of $d_2$. For the pull back bundle $F^*\xi$ over $X_1$, the Thom %space $T(F^*\xi)=T(Y)\cup_{\phi_{F^*\xi}} D^{n+l_1+k}$. Note that the restriction bundle of $F^*\xi$ to $Y$ is %just $\xi_Y$. 
By Lemma \ref{attachThomspace}, one has $c\phi_{\xi}\simeq \Sigma(\Sigma^k d_1)$.
According to Example 2 and Proposition 4 in \cite[Page 24]{frank1967reducible}, we have
\[
[\phi_{F^*\xi}]\in -i_* \{J_{n-1,k}(f_{\xi_Y}),\Sigma ^k d_1,\Sigma^{k-1}d_2 \}.
\] 
\end{proof}

\subsection{Normal bundles of topological manifolds and $\mathrm{PL}$ manifolds}
Recall that for a topological $n$-manifold $M_{\mathrm{Top}}$, there exist a topological embedding $M_{\mathrm{Top}}\hookrightarrow S^{n+k}$ with normal $(\R^k,0)$-bundle $\eta_{M_{\mathrm{Top}}}$, unique up to isotopy, for $k>>n$ (e.g. \cite[Page 203]{kirby1977foundational} \& \cite{kister1964microbundles}). Since $\pi_i B\T_0(D^k)\cong \pi_i B\T_0 \R^k$ when $k>>i$, there exists a $(D^k,0)$-bundle,  unique up to isotopy, $\eta^{D^k}_{M_{\mathrm{Top}}}$ whose interior points are isomorphic to the topological $(\R^k,0)$ bundle $\eta_{M_{\mathrm{Top}}}$ when $k$ is very large. 

Note that the Thom space $T(\eta^{D^k}_{M_{\mathrm{Top}}})$ is homeomorphic to the Pontrjagin-Thom construction
\[
S^{n+k}/(S^{n+k}-\eta_{M_{\mathrm{Top}}})
\]
of the topological embedding $M_{\mathrm{Top}}\hookrightarrow S^{n+k}$.
 
For any compact $\mathrm{PL}$ manifold $M_{\mathrm{PL}}$, there also exists a unique $\mathrm{PL}$ $(D^k,0)$-bundle $\eta^{D^k}_{M_{\mathrm{PL}}}$ whose interior points are isomorphic to the $\mathrm{PL}$ normal bundle $\eta_{M_{\mathrm{Top}}}$ when $k$ is very large (\cite[Theorem 4 \& Corollary 4.2]{haefliger1965piecewise}). 
Note that $\mathrm{PL}$ $(D^k,0)$-bundles are also topological $(D^k,0)$-bundles through the forgetful map
\[
[M_{\mathrm{PL}},B\mathrm{PL}_0(D^k)]\longrightarrow [M_{\mathrm{PL}},B\mathrm{Top}_0(D^k)].
\]
The Thom space $T(\eta^{D^k}_{M_{\mathrm{PL}}})$ can be regarded as the Thom space of its forgetful topological $(D^k,0)$-bundle.
\subsection{Frank's invariant set $\Delta(M^n)$}
Let $M^n$ be a compact, connected, $n$-dimensional smooth manifold whose boundary $\partial M^n$ is a homotopy $(n-1)$-sphere. By attaching a cone $\mathrm{Cone}(\partial M^n)$ on the boundary of $M^n$, one obtains a closed $\mathrm{PL}$-manifold $\hat{M}^n=M^n\cup\mathrm{Cone}(\partial M)$.
Let $\eta_{M^n}$ be the stable normal vector bundle of $M^n$. 
\begin{lem}[Frank, {\cite[Corollary 1, Page 19]{frank1967reducible}}]\label{Franklem1exist}
	There exists a vector bundle $\xi_{\hat{M}^n}$ over $\hat{M}^n$ such that the restriction of $\xi_{\hat{M}^n}$ to $M$ is stably isomorphic to $\eta_{M^n}$.
\end{lem}
\begin{proof}
See the proof in \cite[Corollary 1, Page 19]{frank1967reducible}. 
\end{proof}
%Choose such a vector bundle $\eta_{\hat{M}^n}$ over $\hat{M}^n$ whose restriction to $M^n$ is stably isomorphic to $\eta_{M^n}$.
Suppose the dimension of the stable normal $\mathrm{PL}$-bundle $\eta_{\hat{M}^n}$ is $k>>n$. The Pontrjagin-Thom construction of the $\mathrm{PL}$-embedding implies the Thom space $T(\eta_{\hat{M}^n})$ is reducible (\cite[Page 12]{frank1967reducible}). Let $\xi_{\hat{M}^n}$ be a vector bundle over $\hat{M}^n$ with $\xi_{\hat{M}^n}|_{M^n}\cong \eta_{M^n}$. In \cite[Section 3, Page 563]{frank1968invariant}, Frank defined a subset $\Delta'(\xi_{\hat{M}^n})$ of $\pi_{k+n-1}S^{k}$ to be
\[
\Delta_k'(\xi_{\hat{M}^n})\colon =i_{\ast}^{-1}[\phi_{\xi_{\hat{M}^n}}]=\{x\in \pi_{k+n-1}S^{k} \ | \ i_{\ast}x= [\phi_{\xi_{\hat{M}^n}}]\}.
\]

\begin{lem}[Frank, {\cite[Remark, Page 20]{frank1967reducible}}]\label{FranklemImJ}
	Let $\xi_{\alpha}$ and $\xi_{\beta}$ be two vector bundles over $\hat{M}^n$ with $\xi_{\alpha}|_{M^n}\cong \xi_{\beta}|_{M^n}\cong \eta_{M^n}$, then
	\begin{itemize}
		\item[(a)]$
		\Delta_k'(\xi_{\alpha})+\mathrm{Im}J_{n-1,k}=\Delta'(\xi_{\beta})+\mathrm{Im}J_{n-1,k}$
		\item[(b)] $\Sigma	(\Delta_k'(\xi_{\alpha}))=\Delta_{k+1}'(\xi_{\alpha}\oplus \underline{\R})$.
	\end{itemize}

\end{lem}
\begin{proof}
By Frank's formula \cite[Proposition 1, Page 17]{frank1967reducible}, $[\phi_{\alpha}]-[\phi_{\beta}]=i_{\ast}J^{\mathrm{T}}_{n-1}\Theta_{(\xi_{\alpha},\xi_{\beta})}$. 
When $\xi_{\alpha}$ and $\xi_{\beta}$ are vector bundles, $\Theta_{(\xi_{\alpha},\xi_{\beta})}|_{x\times D^{k}}\in O(k)$ for $x\in D^n$ and $J^{\mathrm{T}}_{n-1}\Theta_{(\xi_{\alpha},\xi_{\beta})}$ lies in the image of $J_{n-1,k}$. Hence, we have the equality $\Delta'(\xi_{\alpha})+\mathrm{Im}J_{n-1,k}=\Delta'(\xi_{\beta})+\mathrm{Im}J_{n-1,k}$.

For the bundle $\xi_{\alpha}\oplus \underline{\R}$, the Thom space $T(\xi_{\alpha}\oplus \underline{\R})\cong \Sigma T (\xi_{\alpha})\cong  T(\eta_{M^n}\oplus  \underline{\R})\cup_{\Sigma \phi_{\alpha}} D^{k+1}$ and one has
\[
\begin{CD}
	\pi_{n-1+k} S^{k}  @>i_*>>  \pi_{n-1+k} T(\xi_{\alpha})\\
	@V\Sigma VV         @V  \Sigma VV\\
		\pi_{n+k} S^{k+1}  @>i_*>>  \pi_{n+k} T(\xi_{\alpha}\oplus \underline{\R}).
\end{CD}
\] When $k$ is very large, $\Sigma$ is bijective and one has $\Sigma	(\Delta_k'(\xi_{\alpha}))=\Delta_{k+1}'(\xi_{\alpha}\oplus \underline{\R})$.
\end{proof}
Recall that $q$ is the quotient map $\pi_{n-1}^s \longrightarrow \coker J_{n-1}=\pi^s_{n-1}/\mathrm{Im}J_{n-1}$.
\begin{defn}[Frank, {\cite[Section 3, Page 563]{frank1968invariant}}]
Let $M^n$ be a compact smooth manifold of dimension $n$ whose boundary is a homotopy sphere and let $\eta_{M^n}$ be the stable normal vector bundle of $M^n$ of dimension $k$ with $k$ very large.
Define
\[
\Delta(M^n):=q(\Sigma^{\infty}\Delta_k'(\xi_{\hat{M}^n}))\subset \pi^s_{n-1}/\mathrm{Im}J_{n-1},
\]
where $\xi_{\hat{M}^n}$ is a vector bundle over $\hat{M}^n$ of dimension $k$ with $\xi_{\hat{M}^n}|_{M^n}\cong \eta_{M^n}$.
By Lemma \ref{Franklem1exist}, such $\xi_{\hat{M}^n}$ exists and $\Delta(M^n)$ is well-defined by Lemma \ref{FranklemImJ} and Remark \ref{remmarkfibrepreserving}.
\end{defn}

Let $P$ be the quotient map $\Theta_{n-1}\longrightarrow \Theta_{n-1}/bP_{n}\subset  \pi^s_{n-1}/\mathrm{Im}J_{n-1}$. Frank gave a sketch proof of the following theorem:
\begin{thm}[Frank, {\cite[Theorem 1]{frank1968invariant}}]\label{Frankthm1}
$-P([\partial M^n])\in \Delta(M^n)$.
\end{thm}
\begin{proof}
We follow Brumfiel's method in \cite[Lemma 3.1, Page 299]{BrumfielI}. When $k>>n$, the cofibration sequence
\[
\begin{CD}
	M^n @>>>\hat{M}^n @>c_{M^n}>>  S^{n}
\end{CD}
\]
induces the following exact sequences
\[
\begin{CD}
	[S^{n},B\P_0(D^k)] @>c^{\ast}_{M^n}>>	[\hat{M}^n,B\P_0(D^k)] @>r^{PL}_{M^n}>> [M^n,B\P_0(D^k)]\\
	@V\cong V int V             @V\cong V int V             @V\cong V int V\\
 [S^{n},B\P_0(\R^k)] @>c^{\ast}_{M^n}>>	[\hat{M}^n,B\P_0(\R^k)] @>r^{PL}_{M^n}>> [M^n,B\P_0(\R^k)].\\
\end{CD}
\]
The $\mathrm{PL}$-isomorphism $r^{PL}_{M^n}\eta_{\hat{M}^n}\cong r^{PL}_{M^n}\xi_{\hat{M}^n}$ implies there exists a $\mathrm{PL}$ $(D^k,0)$-bundle $\xi^{D^k}_{S^n}$ over $S^n$ with $c^{\ast}_{M^n} (\xi^{D^k}_{S^n})$ stably isomorphic to the difference bundle $\xi^{D^k}_{\hat{M}^n}-\eta^{D^k}_{\hat{M}^n}$. 
According to \cite[4.4, Page 304]{BrumfielI} and \cite[Page 305]{BrumfielII}, one has commutative diagram 
\[
\begin{CD}
 \pi_n B\P @>\beta>> \Theta_{n-1} @>\cong>>  \pi_{n-1}\mathrm{PL}/O \\
          @V B J_{n-1}^{\P}VV             @V P VV    @VVV\\
\pi^s_{n-1} \cong\pi_n B\Omega^{\infty}S^{\infty} @>q>>  \pi^s_{n-1}/\mathrm{Im}J_{n-1} @>>>  \pi_{n-1}(\Omega^{\infty}S^{\infty})/O
\end{CD}
\]
and Brumfiel‘s method in the proof of Lemma 3.1 in \cite{BrumfielI,BrumfielII} shows that $\beta(int(\xi^{D^k}_{S^n}))=-[\partial M^n]$. 

On the other hand, when $k>>n$ we can decompose $B J_{n-1}^{\P}$ as follows
\[
\begin{CD}
		\pi_n B\P  @>>> \pi_n B\T_0   @>>> \pi_n B(\Omega^{\infty}S^{\infty})\\
	@AAA                 @AA\cong A    @|\\
	\pi_n B\P_0 \R^k  @>>> \pi_n B\T_0\R^k   @>>> \pi_n B(\Omega^{\infty}S^{\infty})\\
	@AAA                 @AA\cong A    @|\\
		\pi_n B\P_0 D^k  @>>> \pi_n B\T_0 D^k   @>>> \pi_n B(\Omega^{\infty}S^{\infty})\\
		@.            @AA\cong A          @|\\
		              @.   \pi_{n-1} \T_0 D^k   @> \Sigma^{\infty}J^{\mathrm{T}}_{n-1}>>  \pi_{n-1} (\Omega^{\infty}S^{\infty})\cong \pi^s_{n-1}.
\end{CD}
\]
As a topological $(D^k,0)$-bundle, the clutching map of $\xi^{D^k}_{S^n}$ is just $\Theta_{\xi^{D^k}_{\hat{M}^n},\eta^{D^k}_{\hat{M}^n}}$. The one-to-one correspondence $[S^{n},\mathrm{B\T_0 D^k}]\cong [S^{n-1},\T_0 D^k]$ implies the bundle $\xi^{D^k}_{S^n}$ corresponds to the map $\Theta_{\xi^{D^k}_{\hat{M}^n},\eta^{D^k}_{\hat{M}^n}}$.
%\[
%\Theta_{\xi_{\hat{M}^n},\eta_{\hat{M}^n}}\colon S^{n-1}\rightarrow \T_0(D^k)
%\]
%\[
%x\mapsto \mathrm{pr}_2 \tilde{\Theta}_{\xi_{\hat{M}^n},\eta_{\hat{M}^n}}(x,-).
%\]
Then one has
\[
-P[\partial M^n]=P\beta(\xi^{D^k}_{S^n})=q( \Sigma^{\infty}J^{\mathrm{T}}_{n-1}\Theta_{\xi_{\hat{M}^n},\eta_{\hat{M}^n}}).
\]
Since $T(\eta^{D^k}_{\hat{M}^n})$ is reducible and by Lemma \ref{FranklemThomSpace},
\[
i_* [J_{n-1}^{\mathrm{T}}\Theta_{\xi_{\hat{M}^n},\eta_{\hat{M}^n}}]=[\phi_{\xi_{\hat{M}^n}}],
\]
one has $-P[\partial M^n]=q( \Sigma^{\infty}(J^{\mathrm{T}}_{n-1}\Theta_{\xi_{\hat{M}^n},\eta_{\hat{M}^n}}))\in q(\Sigma^{\infty}\Delta'_k( \xi_{\hat{M}^n})) =\Delta(M^n)$.
\end{proof}

\subsection{Proof of Theorem \ref{FrankSmith}}
In this subsection, we prove Theorem \ref{FrankSmith}.
%\begin{thm}[The Frank--Kosinski--Smith theorem]
%	For $\alpha\in \pi_i SO(j)$ and $\beta \in \pi_j SO(i))$ with $i,j\geq 1$, let $u=J_{j,i}\beta$ and suppose
%	\begin{itemize}
%		\item[(1)] $J_{i,j}\alpha= vw$ with $v$ finite order in $\pi_{j+t}S^j$, $ 0<t<i$;
%		\item[(2)] 	$s_{\infty}\alpha=0$;
%	where $s_{\infty}\colon \pi_i SO(j)\longrightarrow \pi_i SO$ is the homomorphisms induced by the natural inclusion;
%		\item[(3)] $(\Sigma^{\infty}u) (\Sigma^{\infty}v)=0$;
%		\item[(4)] $(\Sigma^{\infty}u)\circ \pi^s_{i+1}\subset \mathrm{Image}J_{i+j+1}$,
%	\end{itemize}
%	then one has the following inclusion
%	\[
%	(-1)^{(i+1)(j+1)-1}P(\sigma_{i,j}(\alpha,\beta))\in q(<\Sigma^{\infty}u, \Sigma^{\infty}v, \Sigma^{\infty}w>),
%	\]
%	where $<\Sigma^{\infty} u, \Sigma^{\infty}v, \Sigma^{\infty}w>$ denotes the stable Toda bracket of $\Sigma^{\infty}u$, % $\Sigma^{\infty}v$ and $\Sigma^{\infty}w$.
%\end{thm}
\FrankKosinskiSmiththm*
\begin{proof}
Let $\xi_{s_1\alpha}$ and $\xi_{s_1\beta}$ be the vector bundles corresponding to $s_1\alpha\in \pi_iSO(j+1)$ and
$s_1\beta\in \pi_j SO(i+1)$. Let $M(\alpha,\beta)$ be the plumbing of two disc bundles of $\xi_{s_1\alpha}$ and $\xi_{s_1\beta}$.

Note that $M(\alpha,\beta)\simeq S^{i+1}\vee S^{j+1}$. Let $p_2\colon M(\alpha,\beta)\longrightarrow S^{j+1}$
be the collapsing map. Write $\xi_{-s_1\beta}$ for the stable inverse bundle of $\xi_{s_1\beta}$ of dimension $k$. Since the vector bundle $\xi_{s_1 \alpha}$ is stably trivial, the pull back bundle $p_{2}^*\xi_{-s_1\beta}$, denoted by $\eta_{M}$, is isomorphic to the stable normal bundle of $M(\alpha,\beta)$.

Recall that the composition 
$$p_2i_{\partial}\colon S^{i+j+1}\simeq \partial M(\alpha,\beta) \hookrightarrow M(\alpha,\beta)\longrightarrow S^{j+1}$$
is homotopic to $c_{ij}(J_{i,j+1}s_1\alpha)=c_{ij}(\Sigma J_{i,j}\alpha)=c_{ij}(\Sigma v \Sigma w)$ with $c_{ij}=(-1)^{(i+1)(j+1)}$ (e.g. \cite[Proposition 6.2]{KAHN196581})\footnote{The  expression here differs from that in \cite[Proposition 6.2]{KAHN196581} by -1, because the author of \cite{KAHN196581} used the relation $EJ(\alpha)=-J(\alpha')$ in \cite[Page 91]{KAHN196581}.}. One has the following homotopy commutative diagram
\[
\begin{CD}
	S^{i+j+1} @>i_{\partial} >>  M(\alpha,\beta)  @>>>  \hat{M}(\alpha,\beta)=M(\alpha,\beta) \cup_{i_{\partial}} D^{i+j+2}\\
	@|     @V p_2 VV            @V \hat{p}_2 VV\\
	S^{i+j+1}  @>\Sigma v (c_{ij}\Sigma w)>>  S^{j+1}   @>>>   S^{j+1}\cup_{c_{ij}(J_{i,j+1}s_1\alpha )} D^{i+j+2}=S^{j+1}\cup_{\Sigma v (c_{ij}\Sigma w)} D^{i+j+2}.\\
		@VV c_{ij}\Sigma w V     @|            @V F VV\\
	S^{t+j+1}  @>(\Sigma v )>>  S^{j+1}   @>>>  S^{j+1}\cup_{(\Sigma v )} D^{t+j+2}.\\
\end{CD}
\]
The assumption $v$ is of finite order implies $(s_1\beta)  v $ is also an element of finite order with 
$$J_{t+j} ((s_{\infty}\beta) v)=\Sigma^{\infty}u \Sigma^{\infty}v =0.$$
It follows that $(s_{\infty}\beta) v=0$ and the stable vector bundle $\xi_{-s_1\beta}$ can be extended to a stable vector bundle $\hat{\xi}_{-s_1\beta}$ over  $S^{j+1}\cup_{(\Sigma v )} D^{t+j+2}$ when $k$ is very large. 
The pull-back bundle $F^*\hat{\xi}_{-s_1\beta}$ over $S^{j+1}\cup_{\Sigma v (c_{ij}\Sigma w) } D^{i+j+2}$ is also an extension of $\xi_{-s_1\beta}$. 

Let $\xi_{\hat{M}}=\hat{p_2}^* F^*\hat{\xi}_{-s_1\beta}$ be the pull back bundle over  
$ \hat{M}(\alpha,\beta)$ which extends the stable normal bundle of $M(\alpha,\beta)$.
Consider the Thom spaces of $\xi_{\hat{M}}$ and $F^*\hat{\xi}_{-s_1\beta}$, one has the following homotopy commutative diagram
\[
\begin{CD}
	S^{i+j+1+k}  @> \phi_{\xi_{\hat{M}}}>>  T(\eta_M) @>>> T(\xi_{\hat{M}})=T(\eta_M)\cup_{\phi_{\xi_{\hat{M}}}} D^{i+j+2+k}\\
	@|  @VT(p_2)VV @VT(\hat{p}_2)VV\\
		S^{i+j+1+k}  @> \phi_{F^*\hat{\xi}_{-s_1\beta}}>>  T(\xi_{-s_1\beta})  @>>>   T(F^*\hat{\xi}_{-s_1\beta})=T(\xi_{-s_1\beta})\cup_{\phi_{F^*\hat{\xi}_{-s_1\beta}}} D^{i+j+2+k}.\\
 \end{CD}
\]
For the inclusions of $S^k$ in $ T(\eta_M)$ and  $T(\xi_{-s_1\beta})$, one has
\[
\begin{CD}
S^k  @>i_M>>  T(\eta_M) \\
@|   @V T(p_2)VV \\
S^k  @>i_{S^{j+1}}>> T(\xi_{-s_1\beta}).
\end{CD}
\]
Therefore, for $x\in \Delta_k'(\xi_{\hat{M}})\subset \pi_{i+j+1+k}S^{k}$, one has 
$$i_{S^{j+1},*}(x)=[T(p_2) i_{M,*} (x)]=[T(p_2)  \phi_{\xi_{\hat{M}}} ]=[\phi_{F^*\hat{\xi}_{-s_1\beta}}]\in \pi_{i+j+1+k}T(\xi_{-s_1\beta}).$$
According to Lemma \ref{TodaLemma}, write $X$ for $S^{j+1}\cup_{\Sigma v (c_{ij}\Sigma w)} D^{i+j+2}$ and $X_1$ for $S^{j+1}\cup_{\Sigma v } D^{t+j+2}$. Note that the vector bundle $\hat{\xi}_{-s_1\beta}$ over $X_1$ satisfies $F^*(\hat{\xi}_{-s_1\beta}){|_{S^{j+1}}}=\xi_{-s_1\beta}$. We have
\[
i_{S^{j+1},*} (x)=[\phi_{F^*\hat{\xi}_{-s_1\beta}}]\in -i_{S^{j+1},*}\{  -\Sigma^{k-i}u, \Sigma^{k} v,c_{ij}\Sigma^{k}w  \}\subset i_{S^{j+1},*}(c_{ij}\{  \Sigma^{k-i}u, \Sigma^{k} v,\Sigma^{k}w  \}).
\]
On the other hand, by the Blakers--Massey theorem (see e.g. \cite[Proposition 4.28]{hatcher2002algebraic}), when $k$ is very large the cofiber sequence 
\[
S^{j+k} \longrightarrow S^k  \longrightarrow  S^k\cup_{J_{j,k}(-s_1\beta)}D^{j+1+k}\cong T(\xi_{-s_1\beta}) 
\]
induces the following exact sequences of homotopy groups
\[
\footnotesize{
\begin{CD}
	\pi_{i+j+1+k} S^{j+k} @> J_{j,k}(-s_1\beta)>> 	\pi_{i+j+1+k} S^{k}    @>i_{S^{j+1}}>>  \pi_{i+j+1+k}  S^k\cup_{J_{j,k}(-s_1\beta)}D^{j+k+1}\\
	@|        @|      @A \cong AA\\
		\pi_{i+j+1+k} S^{j+k} @>J_{j,k}(-s_1\beta) >> 	\pi_{i+j+1+k} S^{k}    @>>>  \pi_{i+j+1+k} (S^k\cup_{J_{j,k}(-s_1\beta)}S^{j+k}\times [0,1],S^{j+k}).
\end{CD}}
\]
We see $\mathrm{Ker}( i_{S^{j+1}})\subset  J_{j,k}(-s_1\beta)\circ 	\pi_{i+j+1+k} S^{j+k} \subset  (\Sigma^{k-i}u)\circ 	\pi_{i+j+1+k} S^{j+k}$, and we have  
\[
\Delta_k'(\xi_{\hat{M}})\subset  c_{ij}\{  \Sigma^{k-i}u, \Sigma^{k} v,\Sigma^{k} w \}+  (\Sigma^{k-i}u)\circ 	\pi_{i+j+1+k} S^{j+k}.
\]
By Theorem \ref{Frankthm1}, we have 
\begin{align*}
	-c_{ij}P[\partial M]\in c_{ij}\Delta(M) & = c_{ij} q(\Delta_k'(\xi_{\hat{M}}))\\
	                                            &\subset q(\Sigma^{\infty} \{ \Sigma^{k-i}u, \Sigma^{k} v,\Sigma^{k} w \}+\Sigma^{\infty}u\circ \pi^s_{i+1})\\
	                                            & \subset q(<\Sigma^{\infty}u, \Sigma^{\infty} v,\Sigma^{\infty} w>).
\end{align*}

\end{proof}

\subsection{Further remarks}
\begin{rem}
In this paper, the proof of Theorem \ref{FrankSmith} relies heavily on classical homotopy theory and classical surgery theory. According to \cite{kosinski1971toda} and \cite[Section 5]{smith1974framings}, there should exist a more geometric proof using Kosinski's approach in \cite{kosinski1971toda}. The key point is to remove the dimension restrictions in Theorem 4.3 (b) of \cite{kosinski1971toda}.
\end{rem}
\begin{rem}\label{remarkcollarcondition}
	In fact, in our paper $\bD(D^n):=|\bD(D^n)^{Collared}_{\bullet}|$ is the  geometric realization of the semi-simplicial group $\bD(D^n)^{Collared}_{\bullet}$.
	Note that in \cite{hebestreit2021vanishing,krannich2022homological}, the definition of the semi-simplicial group $\bD(M)^{Collared}_{\bullet}$ should satisfy the collared condition (\cite[\S 1.3]{krannich2022homological}) to ensure $\bD(M)_{\bullet}^{Collared}$ is Kan \cite[Remark 2.2.2]{hebestreit2021vanishing}. Under this collared condition, many classical results and tools used in this paper still work. We refer to \cite{goodwillie2025functoriality} for more discussion on this topic.
	%In this paper, we won't distinguish $\D(M)$ and $|\D (M)^{Collared}_{\bullet}|$ for convenience.
\end{rem}

\appendix

\section{The Gromoll filtration: table of values}\label{appdixA}
As a supplement to the tables in \cite[Appendix A]{crowley2013gromoll} and \cite[Appendix A]{crowley2018harmonic}, we update some Gromoll filtration groups in the following table.

\begin{longtable}{lcc}
		\toprule
		$n+1$ & $\Gamma^{n+1}_{2}\cong \Theta_{n+1}$ & $\Gamma_{i+1}^{n+1}$ for $i \geq 2$ \\
		\midrule
		%		7  & $\Z_{28}$  &   $\Gamma_{3}^{7}\subset\Z_{14}\subset \Gamma_{2}^{7}$     \\
		%	8  & $\Z_{2}$  &  nothing known\\ 
		9  & $(\Z_{2})^3$  & $(\Z_{2})^2\subset \Gamma_3^9 $ \ by Theorem \ref{Theorem A} (b) \\
		\midrule
		10  & $\Z_{2}\oplus \Z_3$  & $\Gamma_4^{10}=\Gamma_{2}^{10}$ \ by Theorem \ref{Theorem A} (a) \\ 
			\midrule
		%	11  & $\Z_{992}$  &   $\Z_4\subset\Gamma_{5}^{11}\subset \Gamma_{3}^{11}\subset \Z_{496}$     \\
		%	12  & $0$  &  nothing need to know\\ 
		13  & $\Z_{3}$  & $\Gamma_{7}^{13}=\Gamma_{2}^{13}$ \ by Theorem \ref{Theorem A} (a)\\ 
			\midrule
		%	14  & $\Z_{2}$  &  nothing known\\
		15  & $\Z_{2}\oplus \Z_{8128}$ & \makecell{$\Z_{2}\oplus 0\subset \Gamma_9^{15}\subset\Gamma_{4}^{15}= \Gamma_{3}^{15} \cong\Z_{2}\oplus \Z_{4064}$ \\ by Theorem \ref{Theorem A} (b), (c) and by \cite{antonelli1970gromoll,weiss1986spheres}}  \\  
		%	16  & $\Z_{2}$  &  nothing known\\
			\midrule
		17  & $(\Z_{2})^4$  & $(\Z_2)^3\subset\Gamma_{5}^{17}$ \ by Theorem \ref{Theorem A} (b)\\
			\midrule
		18  & $\Z_8\oplus \Z_{2}$  & $\Gamma_{5}^{18}=\Gamma_{2}^{18}$ \ by Theorem \ref{Theorem A} (a)\\
			\midrule
		$4k+3\geq 19$  &    $\Theta_{4k+3}$  &  \makecell{$\Gamma_{4}^{4k+3}=\Gamma_{3}^{4k+3}, \Theta_{4k+3}/\Gamma_{4}^{4k+3}\cong\Z_{2}$\\ by Theorem \ref{Theorem A} (c)}   \\
		\bottomrule
			\caption{Some Gromoll filtration groups}
		\label{tableGromoll}
		\end{longtable}

\section{Some information of $\pi_{1}\D(D^n)$}\label{appdixB}
As a supplement to Theorem \ref{thmpi1}, we list some information of $\pi_1\D(D^n)$ in Table \ref{tableDiff}, including some subgroups and splitness of $\pi_1\D(D^n)$ when $6\leq n \leq 15$ and $n=4k\geq 16$.
\begin{longtable}{lcc}
				\toprule
			$n$ & $\pi_1 \D(D^n)$ & Splitness   \\
			\midrule		
			6 & \makecell{ 
			 	$\Z_{2}\theta_1^6\subset \pi_1 \D(D^6) $
				\\
				\\  by Proposition \ref{dim8z2} (c) }& 
			$\xymatrix{ 
					\Z_{2}\theta_1^6 \ar[d]\ar[dr]^{\cong} &  \\
					\pi_1 \D(D^6)   \ar[r] &   \Theta_8
			}$     	\\   	\midrule                        
			7 &  \makecell{ $\pi_2(\frac{\widetilde{\mathrm{Diff}}_{\partial}(D^7)}{\D(D^7)} )\oplus \Z_{2}\kappa_{1}^7 \oplus \Z_{2} \theta_1^7\oplus\Z_{2}\tau_1^7 $\\ 
				\\ by Corollary \ref{CorollaryC} (b) (see \ref{ProofofCoroC})    }       
			&     	${\xymatrix{ 
					\Z_{2}\kappa_{1}^7 \oplus \Z_{2} \theta_1^7\oplus\Z_{2}\tau_1^7\ar[d]\ar[dr]^{\cong} &  & \\
					\pi_1 \D(D^7)   \ar[r] &   \Theta_9	}}$                       \\	\midrule
			8 &   \makecell{ 
				$\Z_{2}\tau_1^8 \oplus\Z_{3}\theta_1^8 \subset \pi_1 \D (D^8)$    \\
		\\	by Theorem \ref{thmpi1} (a) (see \ref{ProofthmB}) }    &  	${ \xymatrix{ 
					\Z_{2}\tau_1^8 \oplus\Z_{3}\theta_1^8 \ar[d]\ar[dr]^{\cong} &  \\
					\pi_1 \D(D^8)   \ar[r] &   \Theta_{10}}}$  \\	\midrule
			9 & \makecell{  $0\rightarrow \Z_2 \rightarrow \pi_1 \D (D^9) \rightarrow \Theta_{11}\rightarrow 0$\\
			\\	by Corollary \ref{CorollaryC} (see \ref{ProofofCoroC})}
		             	&                   \\	\midrule
			10  &    $\Z_2$ \  by \cite[Example 2.3]{wang2024short}    &                     \\	\midrule
			11  &  \makecell{$\Z_{3} \theta_1^{11}\subset \pi_1 \D (D^{11}) $ \\
				\\ 	by Theorem \ref{thmpi1} (a) (see \ref{ProofthmB})  }  &  	${\xymatrix{ 
					\Z_{3} \theta_1^{11}\ar[d]\ar[dr]^{\cong} &  \\
					\pi_1 \D(D^{11})   \ar[r] &   \Theta_{13}}}$                           \\ 	\midrule
			12 &       $\Z_2$ \ by \cite[Theorem 1.1]{wang2024short}   &    \\  	\midrule
			13 &    \makecell{ $\Z_{2}\theta_{1}^{13}\subset 
				   	\pi_1 \D (D^{13}) $   \\
				   	\\ 	by Theorem \ref{thmpi1} (b) (see \ref{proofThmBb}) }  &      	${\xymatrix{ 
				   			\Z_{2} \theta_1^{13}\ar[d]\ar[dr]^{\cong} &  \\
				   			\pi_1 \D(D^{13})   \ar[r] &   \Theta_{15}/bP_{16}}}$                                                                        \\ 	\midrule
			     	14 &   \makecell{  $\Z_{2}\theta_{1}^{14}\subset 
			   \pi_1 \D (D^{14}) $ \\
			   \\ by Proposition \ref{dim14z2} (c) }  &      	${\xymatrix{ 
			   		\Z_{2} \theta_1^{14}\ar[d]\ar[dr]^{\cong} &  \\
			   		\pi_1 \D(D^{14})   \ar[r] &   \Theta_{16}}}$   \\ 	\midrule
		   		    	15 &  \makecell{  $\Z_{2}  \oplus \Z_{2}\kappa_1^{15} \oplus \Z_{2} \tau_1^{15}\oplus \Z_{2}\eta_{1}^{15} \oplus \Z_{2} \theta_{1}^{15} $  \\
		   		    		\\ by Corollary \ref{CorollaryC} (a) (see \ref{ProofofCoroC})}  &            	${\xymatrix{ 
		   				\Z_{2}\kappa_1^{15} \oplus \Z_{2} \tau_1^{15}\oplus \Z_{2}\theta_{1,1}^{15} \oplus \Z_{2} \theta_{1,3}^{15}\ar[d]\ar[dr]^{\cong} &  \\
		   				\pi_1 \D(D^{15})   \ar[r] &   \Theta_{17}}}$   \\  	\midrule
	   				 	$4k\geq 12$ &   $\Theta_{4k+2}$ \ by \cite[Theorem 1.1]{wang2024short}  &    \\                      
			\bottomrule
	\caption{Some information of $\pi_1\D(D^n)$}
\label{tableDiff}
\end{longtable}

\section{Recollections of some data on homotopy groups and exotic spheres}\label{Appendixhomotopyinformation}
For convenience, we recollect some data on homotopy groups and groups of homotopy spheres, which are used in Section \ref{SectionMilnor}.

\subsection{Data on some $\pi_{k+n}S^n$ and $\pi^s_k$}
For the unstable homotopy groups $\pi_{n+k} S^n$, we use the data in Table \ref{tableUnstable}. The $\Z_2$-information of Table \ref{tableUnstable} comes from \cite[Page 39-50,61]{toda1962composition} and the  information of odd primes comes from \cite[Page 178-181,186,189]{toda1962composition}.

In Table \ref{tabbleStable}, we list the data on stable homotopy groups $\pi^s_k$ of spheres together with the image  of stable $J$-homomorphism $\mathrm{Im}J_k$ and the cokernel $\mathrm{Coker}J_k$ for $k=1,3,5$ and $7 \leq k \leq 18$.
The data on $\pi^s_k$ come from \cite[Page 189]{toda1962composition}. Information on the $J$-homomorphism comes from \cite[Page 21-22 ]{ADAMS196621} and \cite[Theorem 1.1.14]{ravenel2003complex}.
We also list some standard notations and relations of the generators (\cite{toda1962composition}) in Table \ref{tabbleStable} as follows
\begin{notation}
	\begin{itemize}
		\item[(1)]  $\eta=\Sigma^{\infty}\eta_{n} $, $\alpha_{1}=\Sigma^{\infty}  \alpha_{1,S^n}$;
		\item[(2)] $	\nu=\Sigma^{\infty} \nu_{n}, 2\nu=\Sigma^{\infty}\nu'$;
		\item[(3)] $\epsilon=\Sigma^{\infty}\epsilon_3 \in <\nu, 2\nu,\eta>$, $\beta_1\in <\alpha_{1},\alpha_{1},\alpha_{1} >$.
		\end{itemize}
\end{notation}
\begin{relation}
\begin{itemize}
	\item[(1)] $\eta \overline{\nu}=\nu^3$, $\eta\sigma=\overline{\nu}+\epsilon$, $\sigma\epsilon=\sigma\overline{\nu}=\sigma\xi=0$, $\sigma\mu=\eta\rho$;
	\item[(2)] $\eta\kappa\in <\nu,2\nu,\epsilon>$, $\eta^*\in <\sigma,2\sigma,\eta>$ and $\nu^*\in <\sigma,2\sigma,\nu>$.
\end{itemize}
\end{relation}

%\begin{proof}
%	The Toda bracket $\epsilon\in <\nu, 2\nu,\eta>$ comes from \cite[(6.1)]{toda1962composition} and the Toda bracket $\eta\kappa\in <\nu,2\nu,\epsilon>$ comes from \cite[Theorem 10.10]{toda1962composition}. The  relations and the other Toda brackets come from Toda's book \cite[Page 189-190]{toda1962composition}.  
%\end{proof}

	\begin{longtable}{llccc}
		\toprule
		$k$& $n$ & $\pi_{k+n} S^n$ & Generators   & Relations   \\
		\midrule
		1 & $n\geq 3$ & $\Z_{2}$ & $\eta_{n} $   &         $\Sigma \eta_{n}=\eta_{n+1}$  \\
		\midrule
		3 &  3 &  $\Z_3 \oplus \Z_{4}$    & $\alpha_{1,S^3}$,  $\nu'$    &      \\
		\midrule
		3 &  4 &  $\Z\oplus   \Z_{3}\oplus \Z_{4}$   &   $\nu_4$,  $\alpha_{1,S^4}$, $\Sigma \nu'$&    $\Sigma  \alpha_{1,S^3}= \alpha_{1,S^{4}} $          \\
		\midrule
		3 &  $n\geq 5$  & $\Z_{3}\oplus \Z_{8}$    &   $
		\alpha_{1,S^n}$, $\nu_n$   &   $  \begin{array}{c}
			\Sigma \nu_{n}=\nu_{n+1}\\
			\Sigma  \alpha_{1,S^n}= \alpha_{1,S^{n+1}} \\
			2\nu_n=\Sigma^{n-2}\nu'\\		 			 	 
		\end{array}      $   \\
		\midrule
		4 &  3 &  $\Z_2$    &  $\nu'\eta_6$    &          \\
		\midrule
		7  & 7   &  $\Z_{120}$    &  $\Z_{8}\sigma'$, $\Z_{3}\alpha_{2,S^{7}}$, $\Z_{5}\alpha_{1,S^{7}}$     &                   \\
		\midrule
		7  & 8   &  $\Z\oplus \Z_{120}$    &  $\Z\sigma_8\oplus\Z_{8} \Sigma\sigma'$, $\Z_{3}\alpha_{2,S^{7}}$, $\Z_{5}\alpha_{1,S^{7}}$     &         $\Sigma^{\infty}\sigma_8=\sigma, \Sigma^2\sigma'=2\Sigma\sigma_8$          \\
		\midrule
		8  & 3   &  $\Z_2$    &  $\epsilon_3$    &             \\
		%		 &    &      &      &             &    \\ 	  
		\bottomrule
	\caption{Some unstable homotopy groups of $S^n$}
	\label{tableUnstable}
\end{longtable}

	\begin{longtable}{lccc}
		\toprule
		$k$& $\pi^s_{k} $  & $\mathrm{Im}J_k$ & $\mathrm{Coker}J_k$ \\
		\midrule
		1  & $\Z_{2} \eta $   &  $\Z_{2}\eta$  &  0       \\
		\midrule
		3  &  $\Z_3 \alpha_{1} \oplus \Z_{8}\nu$     &  $\Z_{3}\alpha_{1}\oplus\Z_{8}\nu$  & 0   \\
		\midrule
		5  &  0    &  0  & 0   \\
		\midrule
		7  & $\Z_{16}\sigma \oplus  \Z_{3}\alpha_{2} \oplus \Z_{5} \alpha_{1,5}$      
		&$\Z_{16}\sigma \oplus  \Z_{3}\alpha_{2} \oplus \Z_{5} \alpha_{1,5}$  &  0 \\
		\midrule
		8  & $\Z_{2}\overline{\nu} \oplus  \Z_{2} \epsilon $      
		& $\Z_{2}(\overline{\nu}+\epsilon)$ &  $\Z_{2} [\epsilon]$ \\
		\midrule
		9 & $\Z_{2}\eta\overline{\nu} \oplus \Z_{2}\eta \epsilon \oplus \Z_{2}\mu$     & $\Z_{2}(\eta\overline{\nu}+\eta \epsilon)$   & $\Z_{2} [ \eta \epsilon ] \oplus  \Z_{2}[\mu]$     \\
		\midrule
		10  & $\Z_{2}\eta\mu \oplus  \Z_{3} \beta_1 $      
		&0 & $\Z_{2}[\eta\mu] \oplus  \Z_{3} [\beta_1]$ \\
		\midrule
		11  & $\Z_{8}\xi \oplus  \Z_{9} \alpha' \oplus \Z_{7}\alpha_{1,7}$      
		&$\Z_{8}\xi \oplus  \Z_{9} \alpha' \oplus \Z_{7}\alpha_{1,7}$ & 0 \\
		\midrule
		13  & $\Z_{3}\alpha_{1}\beta_1$      
		&0 &$\Z_{3}[\alpha_{1}\beta_1]$   \\
		\midrule
		14 & $\Z_{2}\kappa \oplus \Z_{2}\sigma^2$      
		&0 & $\Z_{2}[\kappa] \oplus \Z_{2}[\sigma^2]$  \\
		\midrule
		15 & $\Z_{32}\rho \oplus \Z_{2}\eta\kappa \oplus \Z_{3}\alpha_4\oplus  \Z_{5} \alpha_{2,5}$      
		& $\Z_{32}\rho \oplus \Z_{3}\alpha_4 \oplus \Z_5\alpha_{2,5}$ & $\Z_{2}[\eta \kappa] $  \\
		\midrule
		16 & $\Z_{2}\eta^*\oplus \Z_{2}\eta\rho$      
		&$\Z_{2}\eta\rho$ & $\Z_{2}[\eta^*]$  \\
		\midrule
		17 & $\Z_{2}\eta\eta^* \oplus \Z_2\nu\kappa \oplus \Z_{2}\eta^2\rho \oplus \Z_{2}\overline{\mu}$      
		&$\Z_2\eta^2\rho$ & $\Z_{2}[\eta\eta^*] \oplus \Z_2 [\nu\kappa] \oplus \Z_{2}[\overline{\mu}]$  \\
		\midrule
		18 & $\Z_{8}\nu^* \oplus\Z_{2}\eta \overline{\mu}$      
		&0 & $\Z_{8}[\nu^*] \oplus\Z_{2}[\eta \overline{\mu}]$  \\
		\bottomrule
			\caption{Some $\pi^s_{k}$}
		\label{tabbleStable}
	\end{longtable}

\subsection{Data on groups of homotopy spheres}
For the groups of homotopy spheres $\Theta_n$, we list the data in Table \ref{tableHomotopySpheres} for $8\leq n\leq 18$. These data come from Table \ref{tabbleStable}, \cite[Chapter 12]{luck2024surgery} and \cite[Table 1]{isaksen2023stable}.  Note that we used Brumfiel's splitting theorem (\cite{BrumfielI,BrumfielII,BrumfielIII}, \cite[Theorem 12.70]{luck2024surgery}) in dimension 9, 11, 15, and 17. In these dimensions, we have the algebraic isomorphism $\Theta_{n} \cong \mathrm{Coker}J_n \oplus bP_{n+1}$, and the generators we listed in this table are not unique and not intrinsic in these dimensions. These generators can be any homotopy spheres that maps to the same element in $\mathrm{Coker}J_n$. 

	\begin{longtable}{lccc}
		\toprule
		$n$& $\Theta_{n} $ & $bP_{n+1}$   & $\Theta_{n}/bP_{n+1}\subseteq \mathrm{Coker} J_n$   \\
		\midrule

		8  & $\Z_{2}\Sigma^8_{[\epsilon]}$    &  0    
		& $\Z_2[\epsilon]$ \\
		\midrule
		9	&   $\Z_{2}\Sigma^9_{[ \eta \epsilon ]}\oplus \Z_{2}\Sigma^9_{[\mu]}\oplus \Z_2 \Sigma^9_K$    &  $\Z_{2} \Sigma^9_K$    &    $\Z_{2} [ \eta \epsilon ] \oplus  \Z_{2}[\mu]$                  \\
		\midrule
		10	&   $\Z_{2} \Sigma^{10}_{[\eta\mu] }\oplus \Z_{3}\Sigma^{10}_{[\beta_1] }$    &  0    &      $\Z_{2}[\eta\mu] \oplus \Z_{3}[\beta_1]$   \\
		\midrule
		11	&   $\Z_{992}\Sigma^{11}_{M}$    &   $\Z_{992}\Sigma^{11}_{M}$   &    0   \\
		\midrule
		13	&   $\Z_{3}\Sigma^{13}_{[\alpha_{1}\beta_1]}$    &   0  &   $\Z_3 [\alpha_{1}\beta_1] $  \\
		\midrule
		14	&   $\Z_{2}\Sigma^{14}_{[\kappa]}$    &   0  &   $\Z_2 [\kappa] $  \\
		\midrule
		15	&   $\Z_{8128}\Sigma^{15}_{M}\oplus \Z_2 \Sigma^{15}_{[\eta\kappa]}$    &   $\Z_{8128}\Sigma^{15}_{M}$   &    $\Z_{2}[ \eta\kappa]$   \\
		\midrule
		16  & $\Z_{2}\Sigma^{16}_{[\eta^*]}$    &  0    
		& $\Z_2[\eta^*]$ \\
		\midrule
		17  &$\Z_{2}\Sigma^{17}_{[\eta\eta^*]}\oplus \Z_2 \Sigma^{17}_{[\nu\kappa]} \oplus \Z_{2}\Sigma^{17}_{[\overline{\mu}]} \oplus \Z_{2} \Sigma^{17}_{K}$  &  $\Z_{2} \Sigma^{17}_{K}$    
		& $\Z_{2}[\eta\eta^*] \oplus \Z_2 [\nu\kappa] \oplus \Z_{2}[\overline{\mu}]$ \\
		\midrule
		18  & $\Z_{8}\Sigma^{18}_{[\nu^*]}\oplus \Z_{2} \Sigma^{18}_{[\eta \overline{\mu}]}$    &  0    
		& $\Z_{8}[\nu^*] \oplus\Z_{2}[\eta \overline{\mu}]$ \\
		\bottomrule
	\caption{$\Theta_n$ for  $8\leq n\leq 18$}
	\label{tableHomotopySpheres}
\end{longtable}

\subsection{Data on some $\pi_k SO(n)$ and $\pi_k SO$}
%Let $s^m_*\colon \pi_i SO(n)\longrightarrow \pi_i SO(n+m)$ and .
Some data on $\pi_{k}SO(n)$ and $\pi_k SO$ are listed in Table \ref{tableSO} (see e.g. \cite[Page 106]{husemoller1994fibre} and \cite[Page 161,162]{kervaire1960some}).

	\begin{longtable}{llccc}
		\toprule
		$k$& $n$ & $\pi_{k} SO(n)$ & $\pi_k SO$ & Generators \\
		\midrule
		3 & 3  & $\Z$ &       &           $t_3$  \\
		3&   4 &  $\Z\oplus \Z$    &          &    $s_1 t_3$,    $\gamma_{3,4}$  \\
		3	&    &      &    $\Z$    &        $\gamma_3=s^{\infty}_*\gamma_{3,4}$     \\
		7 & 7  & $\Z$ &       &           $t_7$  \\
		7&   8 &  $\Z\oplus \Z$    &          &    $s_1 t_7$,    $\gamma_{7,8}$  \\
		7	&    &      &    $\Z$    &        $\gamma_7=s^{\infty}_*\gamma_{7,8}$    \\
		\bottomrule
			\caption{Some $\pi_{k}SO(n)$ and $\pi_k SO$}
		\label{tableSO}
	\end{longtable}

We also need a lemma of unstable $J$-homomorphism.

\begin{lem}\label{unstableJ}
	For a prime $p$, let $Z_{(p)}:=\{\frac{a}{b} \ | \ a,b\in \Z, p \nmid b \}$ be the $p$-localization of $\Z$. For the unstable $J$-homomorphism $J_{i,n}\colon \pi_{i} SO(n)\longrightarrow \pi_{i+n} S^{n}$ and its $p$-localization $(J_{i,n})_{(p)}\colon \pi_{i} SO(n)\otimes \Z_{(p)}\longrightarrow \pi_{i+n} S^{n}\otimes \Z_{(p)}$, one has
	\begin{itemize}
		\item[(a)] $J_{3,3}\colon \pi_3 SO(3)\longrightarrow \pi_6 S^3$ is surjective;
		\item[(b)] $(J_{7,7})_{(2)}\colon \pi_7 SO(7)\otimes \Z_{(2)}\longrightarrow \pi_{14} S^7\otimes \Z_{(2)}$ is surjective.
	\end{itemize} 
\end{lem}
\begin{proof}
	By \cite[Theorem 4.15]{hilton1953note}, we see $(J_{3,3})_{(2)}$ and $(J_{7,7})_{(2)}$ are surjective. For $(J_{3,3})_{(3)}$, we refer to \cite[Proposition E]{krishnarao1967unstable}.
\end{proof}

\bibliographystyle{plainnat}
\bibliography{refgromollv3.bib}

\begin{bibdiv}
\begin{biblist}

\bib{ADAMS196621}{article}{
      author={Adams, J.~F.},
       title={On the groups {$J(X)$-$IV$}},
        date={1966},
        ISSN={0040-9383},
     journal={Topology},
      volume={5},
      number={1},
       pages={21\ndash 71},
  url={https://www.sciencedirect.com/science/article/pii/0040938366900048},
}

\bib{antonelli1970gromoll}{article}{
      author={Antonelli, P.~L.},
      author={Burghelea, D.},
      author={Kahn, P.~J.},
       title={Gromoll groups, {$\mathrm{Diff}S^n$} and bilinear constructions
  of exotic spheres},
        date={1970},
     journal={Bull. Amer. Math. Soc.},
      volume={76},
      number={4},
       pages={772\ndash 777},
}

\bib{ANTONELLI1972}{article}{
      author={Antonelli, P.~L.},
      author={Burghelea, D.},
      author={Kahn, P.~J.},
       title={The non-finite homotopy type of some diffeomorphism groups},
        date={1972},
        ISSN={0040-9383},
     journal={Topology},
      volume={11},
      number={1},
       pages={1\ndash 49},
  url={https://www.sciencedirect.com/science/article/pii/0040938372900213},
}

\bib{bier2006detecting}{incollection}{
      author={Bier, T.},
      author={Ray, N},
       title={Detecting framed manifolds in the 8 and 16 stems},
        date={2006},
   booktitle={Geometric applications of homotopy theory i},
   publisher={Springer},
     address={Berlin},
       pages={32\ndash 39},
        note={Proceedings, Evanston, March 21--26, 1977},
}

\bib{browder1966open}{article}{
      author={Browder, W.},
       title={Open and closed disk bundles},
        date={1966},
     journal={Ann. of Math.},
      volume={83},
      number={2},
       pages={218\ndash 230},
}

\bib{BrumfielI}{article}{
      author={Brumfiel, G.},
       title={On the homotopy groups of {${BPL}$} and {${ PL/O}$}},
        date={1968},
        ISSN={0003-486X},
     journal={Ann. of Math. (2)},
      volume={88},
       pages={291\ndash 311},
}

\bib{BrumfielII}{article}{
      author={Brumfiel, G.},
       title={On the homotopy groups of {$BPL$} and {$PL/O$}. {II}},
        date={1969},
        ISSN={0040-9383},
     journal={Topology},
      volume={8},
       pages={305\ndash 311},
}

\bib{BrumfielIII}{article}{
      author={Brumfiel, G.},
       title={The homotopy groups of {$BPL$} and {$PL/O$}. {III}},
        date={1970},
        ISSN={0026-2285},
     journal={Michigan Math. J.},
      volume={17},
       pages={217\ndash 224},
}

\bib{BurgheleaLashof}{article}{
      author={Burghelea, D.},
      author={Lashof, R.},
       title={The homotopy type of the space of diffeomorphisms. {I}, {II}},
        date={1974},
        ISSN={0002-9947},
     journal={Trans. Amer. Math. Soc.},
      volume={196},
       pages={1\ndash 36; ibid. 196 (1974), 37\ndash 50},
}

\bib{BurgheleaLashofRothenberg}{book}{
      author={Burghelea, D.},
      author={Lashof, R.},
      author={Rothenberg, M.},
       title={Groups of {A}utomorphisms of {M}anifolds},
      series={Lecture Notes in Math.},
   publisher={Springer-Verlag, Berlin-New York},
        date={1975},
      volume={473},
        note={With an appendix (``The topological category'') by E. Pedersen},
}

\bib{cerf2006pseudo}{incollection}{
      author={Cerf, J.},
       title={The pseudo-isotopy theorem for simply connected differentiable
  manifolds},
        date={1971},
   booktitle={Manifolds--{A}msterdam 1970: {P}roceedings of the {N}uffic
  {S}ummer {S}chool on {M}anifolds {A}msterdam, {A}ugust 17--29, 1970},
      series={Lecture Notes in Math.},
      volume={197},
   publisher={Springer},
     address={Berlin, Heidelberg},
       pages={76\ndash 82},
        note={Reprinted in 2006},
}

\bib{crowley2013gromoll}{article}{
      author={Crowley, D.},
      author={Schick, T.},
       title={The {G}romoll filtration, {KO--characteristic} classes and
  metrics of positive scalar curvature},
        date={2013},
     journal={Geom. Topol.},
      volume={17},
      number={3},
       pages={1773\ndash 1789},
}

\bib{crowley2018harmonic}{article}{
      author={Crowley, D.},
      author={Schick, T.},
      author={Steimle, W.},
       title={Harmonic spinors and metrics of positive curvature via the
  {G}romoll filtration and {T}oda brackets},
        date={2018},
     journal={J. Topol.},
      volume={11},
      number={4},
       pages={1077\ndash 1099},
}

\bib{crowley2023derivative}{article}{
      author={Crowley, D.},
      author={Schick, T.},
      author={Steimle, W.},
       title={The derivative map for diffeomorphism of disks: an example},
        date={2023},
     journal={Geom. Topol.},
      volume={27},
      number={9},
       pages={3699\ndash 3713},
}

\bib{crowley2017positive}{article}{
      author={Crowley, D.},
      author={Wraith, D.~J.},
       title={Positive {R}icci curvature on highly connected manifolds},
        date={2017},
     journal={J. Differential Geom.},
      volume={106},
      number={2},
       pages={187\ndash 243},
}

\bib{FarrellHsiang1976}{incollection}{
      author={Farrell, F.~T.},
      author={Hsiang, W.~C.},
       title={On the rational homotopy groups of the diffeomorphism groups of
  discs, spheres and aspherical manifolds},
        date={1976},
   booktitle={Algebraic and geometric topology, {Proc. Sympos. Pure Math.}},
      volume={32},
   publisher={Amer. Math. Soc.},
     address={Providence, RI},
}

\bib{frank1967reducible}{book}{
      author={Frank, D.~L.},
       title={Reducible {Thom} complexes and the smoothing problem},
   publisher={Univ. of California, Berkeley},
        date={1967},
}

\bib{frank1968invariant}{article}{
      author={Frank, D.~L.},
       title={An invariant for almost closed manifolds},
        date={1968},
     journal={Bull. Amer. Math. Soc.},
      volume={74},
      number={3},
       pages={562\ndash 567},
}

\bib{goodwillie2025functoriality}{article}{
      author={Goodwillie, T.},
      author={Igusa, K.},
      author={Malkiewich, C.},
      author={Merling, M.},
       title={On the functoriality of the space of equivariant smooth
  $h$-cobordisms},
        date={2025},
     journal={arXiv preprint arXiv:2303.14892v3},
        note={Thorough revision, 91 pages, 21 Apr 2025},
}

\bib{gromoll1966differenzierbare}{article}{
      author={Gromoll, D.},
       title={Differenzierbare strukturen und metriken positiver kr{\"u}mmung
  auf sph{\"a}ren},
        date={1966},
     journal={Math. Ann.},
      volume={164},
      number={4},
       pages={353\ndash 371},
}

\bib{haefliger1965piecewise}{article}{
      author={Haefliger, A.},
      author={Wall, C. T.~C.},
       title={Piecewise linear bundles in the stable range},
        date={1965},
     journal={Topology},
      volume={4},
      number={3},
       pages={209\ndash 214},
}

\bib{hatcher1983proof}{article}{
      author={Hatcher, A.},
       title={A proof of the {Smale} conjecture, {$\mathrm{Diff}S^3 \simeq
  O(4)$}},
        date={1983},
     journal={Ann. of Math.},
      volume={117},
      number={3},
       pages={553\ndash 607},
}

\bib{hatcher2002algebraic}{book}{
      author={Hatcher, A.},
       title={Algebraic {T}opology},
   publisher={Cambridge University Press},
     address={Cambridge},
        date={2002},
        ISBN={9780521791601},
}

\bib{hebestreit2021vanishing}{article}{
      author={Hebestreit, F.},
      author={Land, M.},
      author={L{\"u}ck, W.},
      author={Randal-Williams, O.},
       title={A vanishing theorem for tautological classes of aspherical
  manifolds},
        date={2021},
     journal={Geom. Topol.},
      volume={25},
      number={1},
       pages={47\ndash 110},
}

\bib{hilton1953note}{article}{
      author={Hilton, P.~J.},
      author={Whitehead, J. H.~C.},
       title={Note on the {Whitehead} product},
        date={1953},
     journal={Ann. of Math.},
      volume={58},
      number={3},
       pages={429\ndash 442},
}

\bib{hirsch2012differential}{book}{
      author={Hirsch, M.~W.},
       title={Differential {T}opology},
      series={Graduate Texts in Math.},
   publisher={Springer New York},
        date={2012},
      volume={33},
}

\bib{husemoller1994fibre}{book}{
      author={Husem{\"o}ller, D.},
       title={Fibre {B}undles},
      series={Graduate Texts in Math.},
   publisher={Springer},
        date={1994},
      volume={20},
}

\bib{igusa1988stability}{article}{
      author={Igusa, K.},
       title={The stability theorem for smooth pseudoisotopies},
        date={1988},
     journal={K-Theory},
      volume={2},
      number={1-2},
       pages={1\ndash 355},
}

\bib{isaksen2023stable}{article}{
      author={Isaksen, D.~C.},
      author={Wang, G.},
      author={Xu, Z.},
       title={Stable homotopy groups of spheres: from dimension 0 to 90},
        date={2023},
     journal={Publ. Math. IH{\'E}S},
      volume={137},
      number={1},
       pages={107\ndash 243},
}

\bib{KAHN196581}{article}{
      author={Kahn, P.~J.},
       title={Characteristic numbers and oriented homotopy type},
        date={1965},
        ISSN={0040-9383},
     journal={Topology},
      volume={3},
      number={1},
       pages={81\ndash 95},
  url={https://www.sciencedirect.com/science/article/pii/0040938365900716},
}

\bib{kervaire1960some}{article}{
      author={Kervaire, M.~A.},
       title={Some nonstable homotopy groups of {L}ie groups},
        date={1960},
     journal={Illinois J. Math.},
      volume={4},
      number={2},
       pages={161\ndash 169},
}

\bib{kirby1977foundational}{book}{
      author={Kirby, R.~C.},
      author={Siebenmann, L.~C.},
       title={Foundational essays on topological manifolds, smoothings, and
  triangulations},
   publisher={Princeton Univ. Press},
        date={1977},
      number={88},
}

\bib{kister1964microbundles}{article}{
      author={Kister, J.~M.},
       title={Microbundles are fibre bundles},
        date={1964},
     journal={Ann. of Math.},
      volume={80},
      number={1},
       pages={190\ndash 199},
}

\bib{kosinski1971toda}{article}{
      author={Kosinski, A.},
       title={Toda brackets in differential topology},
        date={1971},
     journal={Comment. Math. Helv.},
      volume={46},
      number={1},
       pages={113\ndash 123},
}

\bib{krannich2022homological}{article}{
      author={Krannich, M.},
       title={A homological approach to pseudoisotopy theory. {I}},
        date={2022},
     journal={Invent. Math.},
      volume={227},
      number={3},
       pages={1093\ndash 1167},
}

\bib{krannich2021diffeomorphisms}{article}{
      author={Krannich, M.},
      author={Randal-Williams, O.},
       title={Diffeomorphisms of discs and the second {W}eiss derivative of
  {BTop(-)}},
        date={2021},
     journal={arXiv preprint arXiv:2109.03500},
}

\bib{krishnarao1967unstable}{article}{
      author={Krishnarao, G.~V.},
       title={Unstable homotopy of {O}(n)},
        date={1967},
     journal={Trans. Amer. Math. Soc.},
      volume={127},
      number={1},
       pages={90\ndash 97},
}

\bib{kupers2025diffeomorphisms}{article}{
      author={Kupers, A.},
      author={Randal-Williams, O.},
       title={On diffeomorphisms of even-dimensional discs},
        date={2025},
     journal={J. Amer. Math. Soc.},
      volume={38},
      number={1},
       pages={63\ndash 178},
}

\bib{Lawson1973remarks}{article}{
      author={Lawson, T.~C.},
       title={Remarks on the pairings of {B}redon, {M}ilnor, and
  {M}ilnor-{M}unkres-{N}ovikov},
        date={1973},
     journal={Indiana Univ. Math. J.},
      volume={22},
      number={9},
       pages={833\ndash 843},
}

\bib{luck2024surgery}{book}{
      author={L{\"u}ck, W.},
      author={Macko, T.},
       title={Surgery {T}heory: {F}oundations},
   publisher={Springer Nature},
        date={2024},
      volume={362},
}

\bib{milnor1959differentiable}{article}{
      author={Milnor, J.~W.},
       title={Differentiable structures on spheres},
        date={1959},
     journal={Amer. J. Math.},
      volume={81},
      number={4},
       pages={962\ndash 972},
}

\bib{randal2023diffeomorphisms}{misc}{
      author={Randal-Williams, O.},
       title={Diffeomorphisms of discs},
      series={ICM—International Congress of Mathematicians},
   publisher={EMS Press, Berlin},
        date={2023},
      volume={4},
}

\bib{ravenel2003complex}{book}{
      author={Ravenel, D.~C.},
       title={Complex cobordism and stable homotopy groups of spheres},
   publisher={Amer. Math. Soc.},
        date={2003},
}

\bib{rognes2002two}{article}{
      author={Rognes, J.},
       title={Two-primary algebraic {K}-theory of pointed spaces},
        date={2002},
     journal={Topology},
      volume={41},
      number={5},
       pages={873\ndash 926},
}

\bib{rognes2003smooth}{article}{
      author={Rognes, J.},
       title={The smooth {W}hitehead spectrum of a point at odd regular
  primes},
        date={2003},
     journal={Geom. Topol.},
      volume={7},
      number={1},
       pages={155\ndash 184},
}

\bib{smale1962structure}{article}{
      author={Smale, S.},
       title={On the structure of manifolds},
        date={1962},
     journal={Amer. J. Math.},
      volume={84},
      number={3},
       pages={387\ndash 399},
}

\bib{smith1974framings}{article}{
      author={Smith, L.},
       title={Framings of sphere bundles over spheres, the plumbing pairing,
  and the framed bordism classes of rank two simple {L}ie groups},
        date={1974},
     journal={Topology},
      volume={13},
      number={4},
       pages={401\ndash 415},
}

\bib{stolz2007hochzusammenhangende}{book}{
      author={Stolz, S.},
       title={Hochzusammenhängende {M}annigfaltigkeiten und ihre {R}änder},
      series={Lecture Notes in Mathematics},
   publisher={Springer},
     address={Berlin},
        date={2007},
      volume={1116},
        ISBN={978-3-540-72924-8},
}

\bib{toda1962composition}{book}{
      author={Toda, H.},
       title={Composition methods in homotopy groups of spheres},
   publisher={Princeton Univ. Press},
        date={1962},
      number={49},
}

\bib{wang2024short}{article}{
      author={Wang, W.},
       title={A short note on {$\pi_1 \mathrm{Diff}_{\partial} (D^{4k})$} for
  {$k\geq 3$}},
        date={2024},
     journal={Proc. Amer. Math. Soc.},
      volume={152},
      number={09},
       pages={4067\ndash 4073},
}

\bib{weiss1986spheres}{article}{
      author={Weiss, M.},
       title={{Sph{\`e}res exotiques et l’espace de {W}hitehead}},
        date={1986},
     journal={C. R. Acad. Sci. Paris S{\'e}r. I Math.},
      volume={303},
      number={19},
       pages={885\ndash 888},
}

\bib{WWsurvey}{incollection}{
      author={Weiss, M.},
      author={Williams, B.},
       title={Automorphisms of manifolds},
        date={2001},
   booktitle={Surveys on surgery theory, {V}ol. 2},
      series={Ann. of Math. Stud.},
      volume={14},
   publisher={Princeton Univ. Press, Princeton, NJ},
       pages={165\ndash 220},
}

\end{biblist}
\end{bibdiv}

\end{document}